\documentclass[a4paper,12 pt]{article}
\usepackage{amsmath,amsthm,amssymb,latexsym,epic}
\usepackage{bezier}
\usepackage[all]{xy}
\usepackage{enumerate}

\newcommand{\tr}{\mathrm{tr}}
\newcommand{\im}{\mathrm{im}}
\newcommand{\IO}{\mathcal{IO}}
\newcommand{\PIP}{\mathcal{PIP}}
\newcommand{\B}{\mathcal{B}}
\newcommand{\I}{\mathcal{I}}
\newcommand{\PB}{\mathcal{PB}}
\newcommand{\IP}{\mathcal{IP}}
\newcommand{\Ch}{\mathcal{CS}}
\newcommand{\GR}{\mathcal{R}}
\newcommand{\GL}{\mathcal{L}}
\newcommand{\GH}{\mathcal{H}}
\newcommand{\GJ}{\mathcal{J}}
\newcommand{\GD}{\mathcal{D}}
\newcommand{\IS}{\mathcal{IS}}
\newcommand{\IOP}{\mathcal{IOP}}
\newcommand{\IT}{\mathcal{IT}}
\newcommand{\Sym}{\mathcal{S}}

\newcommand{\rank}{\mathrm{rank}}
\newcommand{\Aut}{\mathrm{Aut}}
\newcommand{\id}{\mathrm{id}}
\newcommand{\Fix}{\mathrm{Fix}}
\newcommand{\dom}{\mathrm{dom}}
\newcommand{\ran}{\mathrm{ran}}
\newcommand{\m}{\mathrm{min}}
\newcommand{\lmod}{\mid\!}
\newcommand{\rmod}{\!\mid}

\title{On a new approach to the dual symmetric inverse monoid
$\I^{\ast}_X$}
\author{Victor Maltcev}
\date{}
\begin{document}

\maketitle

\begin{abstract}
We construct the \emph{inverse partition semigroup} $\IP_X$,
isomorphic to the \emph{dual symmetric inverse monoid}
$\I^{\ast}_X$, introduced in~\cite{FL}. We give a convenient
geometric illustration for elements of $\IP_X$. We describe all
maximal subsemigroups of $\IP_X$ and find a generating set for
$\IP_X$ when $X$ is finite. We prove that all the automorphisms of
$\IP_X$ are inner. We show how to embed the symmetric inverse
semigroup into the inverse partition one. For finite sets $X$, we
establish that, up to equivalence, there is a unique faithful
effective transitive representation of $\IP_n$, namely to
$\IS_{2^n-2}$. Finally, we construct an interesting
$\GH$-cross-section of $\IP_n$, which is reminiscent of $\IO_n$, the
$\GH$-cross-section of $\IS_n$, constructed in~\cite{Cowan-Reilly}.
\end{abstract}

\newtheorem{theorem}{Theorem}
\newtheorem{lemma}{Lemma}
\newtheorem{proposition}{Proposition}
\newtheorem{corollary}{Corollary}
\newtheorem{definition}{Definition}
\newtheorem{example}{Example}
\newtheorem{remark}{Remark}

\section{Introduction}\label{sec:intro}

The \emph{dual inverse symmetric monoid} $\I^{\ast}_X$ was
introduced in \cite{FL}. It consists of all \emph{biequivalences} on
a set $X$, i.e. all the binary relations $\alpha$ on $X$ that are
both \emph{full}, that is $X\alpha=\alpha X=X$, and
\emph{bifunctional}, that is
$\alpha\circ\alpha^{-1}\circ\alpha=\alpha$. The multiplication in
$\I^{\ast}_X$ is given by:
\begin{equation}\label{eq:FL}
\alpha\beta=\alpha\circ\bigl(\alpha^{-1}\circ\alpha\vee\beta\circ\beta^{-1}\bigr)\circ\beta,
\end{equation}
for $\alpha,\beta\in\I^{\ast}_X$.

In the present paper we introduce the \emph{inverse partition
semigroup} $\IP_X$, isomorphic to $\I^{\ast}_X$ (see
Theorem~\ref{th:Isomorphism}), and investigate some its properties.
The main idea for considering the same semigroup under another point
of view as in \cite{FL} (see definition of $\IP_X$ below) is to
provide a convenient geometric realization for elements of this
semigroup, which will enable us to handle them more easily. Besides,
the semigroup $\IP_X$ naturally arises as an inverse subsemigroup of
the \emph{composition semigroup} $\Ch_X$ (see
Proposition~\ref{pr:IP_X-inverse-subsemigroup}), constructed below,
a generalization of the semigroup $\Ch_n$, introduced
in~\cite{Changchang}. The latter semigroup is close to, so called,
\emph{Brauer-type semigroups}, which were investigated for different
reasons and from different contexts.

The first paper within these investigations, was the work of
Brauer,~\cite{Brauer}, where he introduced the \emph{Brauer
semigroup} $\B_n$ in connection with representations of orthogonal
groups. One more work, where $\B_n$ was studied in connection with
representation theory is \cite{Kerov}. Further work, dedicated to
$\B_n$ are \cite{KMM}, \cite{Maltcev}, \cite{MM}. For example, in
\cite{KMM} all the $\GL$- and $\GR$-cross-sections are described and
in \cite{MM} a presentation for the singular part of $\B_n$ is given
with respect to its minimal generating set. There are several
generalizations of the Brauer semigroup: the \emph{partial Brauer
semigroup} $\PB_n$, introduced in \cite{Mazorchuk-PB}; the
\emph{composition semigroup} $\Ch_n$, appeared in \cite{Changchang};
the \emph{dual symmetric inverse monoid} $\I^{\ast}_X$, introduced
in \cite{FL}; the finite \emph{inverse partition semigroup} $\IP_n$,
appeared in \cite{Maltcev-IP} (which is isomorphic to
$\I^{\ast}_n$); the \emph{partial inverse partition semigroup}
$\PIP_X$, introduced in \cite{KM}. For other papers, dedicated to
these semigroups we refer reader to \cite{F}, \cite{LF},
\cite{Maltcev-PB-C}, \cite{Mazorchuk-End}.

The main purpose of this paper is to investigate some inner
semigroup properties of $\IP_X$, as well as to establish some
connections of $\IP_X$ with other semigroups.

The paper is organized in the following way. In
section~\ref{sec:definitionIP} we define $\IP_X$. After this, in
section~\ref{sec:Isomorphism}, we prove that the constructed
semigroup $\IP_X$ is isomorphic to $\I^{\ast}_X$. In
section~\ref{sec:green} we characterize the Green's relations and
the natural order in $\IP_X$. In section~\ref{sec:some-objects} we
investigate maximal subsemigroups and ideals of $\IP_X$ and define
the \emph{inverse type-preserving semigroup}. In
section~\ref{sec:AutomorphismsIP_X} we describe the automorphism
group $\Aut(\IP_X)$. In section~\ref{sec:Connections} we obtain a
method how to embed the \emph{symmetric inverse semigroup} $\IS_X$
into the inverse partition one. In section~\ref{ef} we obtain that
$\IP_X$ embeds into $\IS_{2^{\lmod X\rmod}-2}$ when $\lmod
X\rmod\in\mathbb{N}\setminus\{1\}$. Finally, in section~\ref{IOP_n}
we define the \emph{inverse ordered partition semigroup} $\IOP_n$,
which behaves similar to the $\GH$-cross-section $\IO_n$ of $\IS_n$,
studied in \cite{Cowan-Reilly}.

Throughout this paper for $S$ a semigroup we denote by $E(S)$ the
set of all idempotents of $S$. The natural order on an inverse
semigroup $S$ will be denoted by $\leq$, i.e., $a\leq b$ for $a,b\in
S$ if and only if there is an idempotent $e$ of $S$ such that $a=be$
(see~\cite{Howie}). We will also need the notion of the \emph{trace}
$\tr(S)$ of an inverse semigroup $S$: the set $S$ together with the
partial multiplication $\ast$, defined as follows: $a\ast b$ is
defined precisely when $ab\in\GR_{a}\cap\GL_{b}$ and is equal then
to $ab$ (see \cite{Meakin} and section XIV.2 of \cite{Petrich}).
Finally, we recall one more definition. For any inverse semigroup
$S$, the \emph{inductive groupoid} of $S$, or \emph{imprint}
$\im(S)$ of $S$, is the triple $\bigl(\tr(S),\leq,\star\bigr)$,
where $\leq$ is the natural partial order in $S$, and $\star$ is a
partial product defined by: for $e\in E(S)$, $a\in S$, $e\leq
aa^{-1}$, $e\star a=ea$ (see section XIV.3.4 of \cite{Petrich}).

\section{Definition of the inverse partition semigroup
$\IP_X$} \label{sec:definitionIP}

Throughout all the paper let $X$ be an arbitrary set. We consider a
map $':X\to X'$ as a fixed bijection and will denote the inverse
bijection by the same symbol, that is $(x')'=x$ for all $x\in X$. We
are going to construct a semigroup $\Ch_X$.

Let $\Ch_X$ be the set of all partitions of $X\cup X'$ into nonempty
blocks. If $X\cup X'=\bigcup\limits_{i\in I}^{\cdot}A_i$ is a
partition of $X\cup X'$ into nonempty blocks $A_i$, $i\in I$,
corresponding to an element $a\in\Ch_X$, then we will write
$a=\bigl(A_i\bigr)_{i\in I}$. In the case when
$I=\{i_1,\ldots,i_k\}$ is finite, we will also write
$a=\bigl\{A_{i_1},\ldots,A_{i_k}\bigr\}$.

For $a\in\Ch_X$ and $x,y\in X\cup X'$, we set $x\equiv_{a}y$
provided that $x$ and $y$ are at the same block of $a$. Clearly, we
can realize $a\in\Ch_X$ as the equivalence relation $\equiv_{a}$.
Thus in spite of  the fact that elements of $\Ch_X$ will be
partitions, we will sometimes treat with them as with the associated
equivalence relations.

Take now $a,b\in\Ch_X$. Define a new equivalence relation, $\equiv$,
on $X\cup X'$ as follows:
\begin{itemize}
\item
for $x,y\in X$ we have $x\equiv y$ if and only if $x\equiv_{a} y$ or
there is a sequence, $c_1,\ldots,c_{2s}$, $s\geq 1$, of elements in
$X$, such that $x\equiv_{a} c'_1$, $c_1\equiv_{b} c_2$,
$c'_2\equiv_{a} c'_3,\ldots,$ $c_{2s-1}\equiv_{b} c_{2s}$, and
$c'_{2s}\equiv_{a} y$;
\item
for $x,y\in X$ we have $x'\equiv y'$ if and only if $x'\equiv_{b}
y'$ or there is a sequence, $c_1,\ldots,c_{2s}$, $s\geq 1$, of
elements in $X$, such that  $x'\equiv_{b} c_1$, $c'_1\equiv_{a}
c'_2$, $c_2\equiv_{b} c_3,\ldots,$ $c'_{2s-1}\equiv_{a} c'_{2s}$,
and $c_{2s}\equiv_{b} y'$;
\item
for $x,y\in X$ we have $x\equiv y'$ if and only if $y'\equiv x$ if
and only if there is a sequence, $c_1,\ldots$, $c_{2s-1}$, $s\geq
1$, of elements in $X$, such that $x\equiv_{a} c'_1$, $c_1\equiv_{b}
c_2$, $c'_2\equiv_{a} c'_3,\ldots,$ $c'_{2s-2}\equiv_{a} c'_{2s-1}$,
and $c_{2s-1}\equiv_{b} y'$.
\end{itemize}

\begin{proposition}\label{pr:well-defined}
$\equiv$ is an equivalence relation on $X\cup X'$.
\end{proposition}

\begin{proof}
It follows immediately from the definition of $\equiv$ that this
relation is reflexive and symmetric. Let now $x\equiv y$ and
$y\equiv z$ for some $x,y,z\in X\cup X'$. We are going to establish
that $x\equiv z$. In the rest of the proof we may assume that $y\in
X$, the other case is treated analogously. We have four possible
cases.

\emph{Case 1}. $x,z\in X$. If $x\equiv_{a}y$ or $y\equiv_{a}z$ then
since $\equiv_a$ is an equivalence relation, we immediately obtain
from the definition of $\equiv$ that $x\equiv z$. Otherwise we have
that there exist $c_1,\ldots,c_{2s},d_1,\ldots,d_{2t}$, elements of
$X$, such that $x\equiv_{a} c'_1$, $c_1\equiv_{b} c_2$,
$c'_2\equiv_{a} c'_3,\ldots,$ $c_{2s-1}\equiv_{b} c_{2s}$,
$c'_{2s}\equiv_{a} y$ and $y\equiv_{a} d'_1$, $d_1\equiv_{b} d_2$,
$d'_2\equiv_{a} d'_3,\ldots,$ $d_{2t-1}\equiv_{b} d_{2t}$,
$d'_{2t}\equiv_{a} z$. Now, using transitiveness of $\equiv_a$, we
can write $c'_{2s}\equiv_a d'_1$ and hence $x\equiv z$.

\emph{Case 2}. $x,z\in X'$. Then there are $c_1,\ldots$, $c_{2s-1}$,
$d_1,\ldots$, $d_{2t-1}$, elements of $X$, such that
$x\equiv_{b}c_{2s-1}$, $c'_{2s-1}\equiv_{a}c'_{2s-2}\ldots,$
$c'_3\equiv_{a}c'_2$, $c_2\equiv_{b}c_1$, $c'_1\equiv_{a}y$ and
$y\equiv_{a} d'_1$, $d_1\equiv_{b} d_2$, $d'_2\equiv_{a}
d'_3,\ldots,$ $d'_{2t-2}\equiv_{a} d'_{2t-1}$, $d_{2t-1}\equiv_{b}
z$. Again, using transitiveness of $\equiv_{a}$, we obtain that
$c'_1\equiv_{a}d'_1$, whence $x\equiv z$.

\emph{Case 3}. $x\in X$ and $z\in X'$. There exist $d_1,\ldots$,
$d_{2t-1}$, elements of $X$, such that $y\equiv_{a} d'_1$,
$d_1\equiv_{b} d_2$, $d'_2\equiv_{a} d'_3,\ldots,$
$d'_{2t-2}\equiv_{a} d'_{2t-1}$, and $d_{2t-1}\equiv_{b} z$. If
$x\equiv_{a}y$ then due to transitiveness of $\equiv_{a}$, we have
$x\equiv z$. Otherwise there are $c_1,\ldots,c_{2s}$, elements of
$X$, such that $x\equiv_{a} c'_1$, $c_1\equiv_{b} c_2$,
$c'_2\equiv_{a} c'_3,\ldots,$ $c_{2s-1}\equiv_{b} c_{2s}$, and
$c'_{2s}\equiv_{a} y$. Then it remains to notice that
$c'_{2s}\equiv_{a}d'_1$.

\emph{Case 4}. $x\in X'$ and $z\in X$. Then, since $z\equiv y$ and
$y\equiv x$, according to Case 3, we have that $z\equiv x$, whence
$x\equiv z$.

The proof is complete.
\end{proof}

Thus $\equiv$ defines a partition of $X\cup X'$ into disjoint blocks
and so belongs to $\Ch_X$. Set this partition to be a product
$a\cdot b$ in $\Ch_X$. One can easily show that
$\bigl(\Ch_X,\cdot\bigr)$ is a semigroup. We will call this
semigroup the \emph{composition semigroup} on the set $X$.

\begin{figure}
\special{em:linewidth 0.4pt} \unitlength 0.80mm
\linethickness{0.4pt}
\begin{picture}(150.00,75.00)
\put(20.00,00.00){\makebox(0,0)[cc]{$\bullet$}}
\put(20.00,10.00){\makebox(0,0)[cc]{$\bullet$}}
\put(20.00,20.00){\makebox(0,0)[cc]{$\bullet$}}
\put(20.00,30.00){\makebox(0,0)[cc]{$\bullet$}}
\put(20.00,40.00){\makebox(0,0)[cc]{$\bullet$}}
\put(20.00,50.00){\makebox(0,0)[cc]{$\bullet$}}
\put(20.00,60.00){\makebox(0,0)[cc]{$\bullet$}}
\put(20.00,70.00){\makebox(0,0)[cc]{$\bullet$}}
\put(45.00,10.00){\makebox(0,0)[cc]{$\bullet$}}
\put(45.00,20.00){\makebox(0,0)[cc]{$\bullet$}}
\put(45.00,30.00){\makebox(0,0)[cc]{$\bullet$}}
\put(45.00,40.00){\makebox(0,0)[cc]{$\bullet$}}
\put(45.00,50.00){\makebox(0,0)[cc]{$\bullet$}}
\put(45.00,60.00){\makebox(0,0)[cc]{$\bullet$}}
\put(45.00,70.00){\makebox(0,0)[cc]{$\bullet$}}
\put(45.00,00.00){\makebox(0,0)[cc]{$\bullet$}}

\put(70.00,00.00){\makebox(0,0)[cc]{$\bullet$}}
\put(70.00,10.00){\makebox(0,0)[cc]{$\bullet$}}
\put(70.00,20.00){\makebox(0,0)[cc]{$\bullet$}}
\put(70.00,30.00){\makebox(0,0)[cc]{$\bullet$}}
\put(70.00,40.00){\makebox(0,0)[cc]{$\bullet$}}
\put(70.00,50.00){\makebox(0,0)[cc]{$\bullet$}}
\put(70.00,60.00){\makebox(0,0)[cc]{$\bullet$}}
\put(70.00,70.00){\makebox(0,0)[cc]{$\bullet$}}
\put(95.00,10.00){\makebox(0,0)[cc]{$\bullet$}}
\put(95.00,20.00){\makebox(0,0)[cc]{$\bullet$}}
\put(95.00,30.00){\makebox(0,0)[cc]{$\bullet$}}
\put(95.00,40.00){\makebox(0,0)[cc]{$\bullet$}}
\put(95.00,50.00){\makebox(0,0)[cc]{$\bullet$}}
\put(95.00,60.00){\makebox(0,0)[cc]{$\bullet$}}
\put(95.00,70.00){\makebox(0,0)[cc]{$\bullet$}}
\put(95.00,00.00){\makebox(0,0)[cc]{$\bullet$}}

\drawline(45.00,00.30)(50.00,00.30)
\drawline(55.00,00.30)(60.00,00.30)
\drawline(65.00,00.30)(70.00,00.30)
\drawline(45.00,10.30)(50.00,10.30)
\drawline(55.00,10.30)(60.00,10.30)
\drawline(65.00,10.30)(70.00,10.30)
\drawline(45.00,20.30)(50.00,20.30)
\drawline(55.00,20.30)(60.00,20.30)
\drawline(65.00,20.30)(70.00,20.30)
\drawline(45.00,30.30)(50.00,30.30)
\drawline(55.00,30.30)(60.00,30.30)
\drawline(65.00,30.30)(70.00,30.30)
\drawline(45.00,40.30)(50.00,40.30)
\drawline(55.00,40.30)(60.00,40.30)
\drawline(65.00,40.30)(70.00,40.30)
\drawline(45.00,50.30)(50.00,50.30)
\drawline(55.00,50.30)(60.00,50.30)
\drawline(65.00,50.30)(70.00,50.30)
\drawline(45.00,60.30)(50.00,60.30)
\drawline(55.00,60.30)(60.00,60.30)
\drawline(65.00,60.30)(70.00,60.30)
\drawline(45.00,70.30)(50.00,70.30)
\drawline(55.00,70.30)(60.00,70.30)
\drawline(65.00,70.30)(70.00,70.30)

\drawline(18.00,-03.00)(18.00,13.00)
\drawline(18.00,13.00)(47.00,03.00)
\drawline(47.00,03.00)(47.00,-03.00)
\drawline(47.00,-03.00)(18.00,-03.00)

\drawline(18.00,18.00)(18.00,32.00)
\drawline(18.00,32.00)(22.00,32.00)
\drawline(22.00,32.00)(22.00,18.00)
\drawline(22.00,18.00)(18.00,18.00)

\drawline(18.00,68.00)(18.00,73.00)
\drawline(18.00,73.00)(22.00,73.00)
\drawline(22.00,73.00)(22.00,68.00)
\drawline(22.00,68.00)(18.00,68.00)

\drawline(18.00,48.00)(18.00,62.00)
\drawline(18.00,62.00)(47.00,73.00)
\drawline(47.00,73.00)(47.00,48.00)
\drawline(47.00,48.00)(18.00,48.00)

\drawline(43.00,42.00)(47.00,42.00)
\drawline(47.00,42.00)(47.00,18.00)
\drawline(47.00,18.00)(43.00,18.00)
\drawline(43.00,18.00)(43.00,42.00)

\drawline(18.00,38.00)(18.00,42.00)
\drawline(18.00,42.00)(22.00,42.00)
\drawline(22.00,42.00)(43.00,12.00)
\drawline(43.00,12.00)(47.00,12.00)
\drawline(47.00,12.00)(47.00,08.00)
\drawline(47.00,08.00)(43.00,08.00)
\drawline(43.00,08.00)(22.00,38.00)
\drawline(22.00,38.00)(18.00,38.00)

\drawline(68.00,-03.00)(68.00,03.00)
\drawline(68.00,03.00)(97.00,03.00)
\drawline(97.00,03.00)(97.00,-03.00)
\drawline(97.00,-03.00)(68.00,-03.00)

\drawline(68.00,08.00)(68.00,32.00)
\drawline(68.00,32.00)(72.00,32.00)
\drawline(72.00,32.00)(72.00,08.00)
\drawline(72.00,08.00)(68.00,08.00)

\drawline(68.00,38.00)(68.00,42.00)
\drawline(68.00,42.00)(97.00,42.00)
\drawline(97.00,42.00)(97.00,38.00)
\drawline(97.00,38.00)(68.00,38.00)

\drawline(68.00,58.00)(68.00,62.00)
\drawline(68.00,62.00)(72.00,62.00)
\drawline(72.00,62.00)(72.00,58.00)
\drawline(72.00,58.00)(68.00,58.00)

\drawline(68.00,67.00)(68.00,73.00)
\drawline(68.00,73.00)(97.00,73.00)
\drawline(97.00,73.00)(97.00,67.00)
\drawline(97.00,67.00)(68.00,67.00)

\drawline(93.00,62.00)(97.00,62.00)
\drawline(97.00,62.00)(97.00,48.00)
\drawline(97.00,48.00)(93.00,48.00)
\drawline(93.00,48.00)(93.00,62.00)

\drawline(93.00,22.00)(97.00,22.00)
\drawline(97.00,22.00)(97.00,08.00)
\drawline(97.00,08.00)(93.00,08.00)
\drawline(93.00,08.00)(93.00,22.00)

\drawline(68.00,48.00)(68.00,54.00)
\drawline(68.00,54.00)(97.00,33.00)
\drawline(97.00,33.00)(97.00,26.00)
\drawline(97.00,26.00)(68.00,48.00)

\drawline(120.00,-03.00)(120.00,13.00)
\drawline(120.00,13.00)(149.00,03.00)
\drawline(149.00,03.00)(149.00,-03.00)
\drawline(149.00,-03.00)(120.00,-03.00)

\drawline(120.00,18.00)(120.00,32.00)
\drawline(120.00,32.00)(124.00,32.00)
\drawline(124.00,32.00)(124.00,18.00)
\drawline(124.00,18.00)(120.00,18.00)

\drawline(120.00,68.00)(120.00,73.00)
\drawline(120.00,73.00)(124.00,73.00)
\drawline(124.00,73.00)(124.00,68.00)
\drawline(124.00,68.00)(120.00,68.00)

\drawline(145.00,62.00)(149.00,62.00)
\drawline(149.00,62.00)(149.00,48.00)
\drawline(149.00,48.00)(145.00,48.00)
\drawline(145.00,48.00)(145.00,62.00)

\drawline(145.00,22.00)(149.00,22.00)
\drawline(149.00,22.00)(149.00,08.00)
\drawline(149.00,08.00)(145.00,08.00)
\drawline(145.00,08.00)(145.00,22.00)

\drawline(120.00,38.00)(120.00,42.00)
\drawline(120.00,42.00)(149.00,42.00)
\drawline(149.00,42.00)(149.00,38.00)
\drawline(149.00,38.00)(120.00,38.00)

\drawline(120.00,48.00)(120.00,63.00)
\drawline(120.00,63.00)(145.00,73.00)
\drawline(145.00,73.00)(149.00,73.00)
\drawline(149.00,73.00)(149.00,67.00)
\drawline(149.00,67.00)(143.00,66.00)
\drawline(143.00,66.00)(143.00,33.00)
\drawline(143.00,33.00)(149.00,32.00)
\drawline(149.00,32.00)(149.00,27.00)
\drawline(149.00,27.00)(145.00,27.00)
\drawline(145.00,27.00)(124.00,48.00)
\drawline(124.00,48.00)(120.00,48.00)

\put(122.00,00.00){\makebox(0,0)[cc]{$\bullet$}}
\put(122.00,10.00){\makebox(0,0)[cc]{$\bullet$}}
\put(122.00,20.00){\makebox(0,0)[cc]{$\bullet$}}
\put(122.00,30.00){\makebox(0,0)[cc]{$\bullet$}}
\put(122.00,40.00){\makebox(0,0)[cc]{$\bullet$}}
\put(122.00,50.00){\makebox(0,0)[cc]{$\bullet$}}
\put(122.00,60.00){\makebox(0,0)[cc]{$\bullet$}}
\put(122.00,70.00){\makebox(0,0)[cc]{$\bullet$}}
\put(147.00,10.00){\makebox(0,0)[cc]{$\bullet$}}
\put(147.00,20.00){\makebox(0,0)[cc]{$\bullet$}}
\put(147.00,30.00){\makebox(0,0)[cc]{$\bullet$}}
\put(147.00,40.00){\makebox(0,0)[cc]{$\bullet$}}
\put(147.00,50.00){\makebox(0,0)[cc]{$\bullet$}}
\put(147.00,60.00){\makebox(0,0)[cc]{$\bullet$}}
\put(147.00,70.00){\makebox(0,0)[cc]{$\bullet$}}
\put(147.00,00.00){\makebox(0,0)[cc]{$\bullet$}}
\put(57.50,35.00){\makebox(0,0)[cc]{$\cdot$}}
\put(109.50,35.00){\makebox(0,0)[cc]{$=$}}

\put(12.00,00.00){\makebox(0,0)[cc]{$1\rightarrow$}}
\put(12.00,10.00){\makebox(0,0)[cc]{$2\rightarrow$}}
\put(12.00,20.00){\makebox(0,0)[cc]{$3\rightarrow$}}
\put(12.00,30.00){\makebox(0,0)[cc]{$4\rightarrow$}}
\put(12.00,40.00){\makebox(0,0)[cc]{$5\rightarrow$}}
\put(12.00,50.00){\makebox(0,0)[cc]{$6\rightarrow$}}
\put(12.00,60.00){\makebox(0,0)[cc]{$7\rightarrow$}}
\put(12.00,70.00){\makebox(0,0)[cc]{$8\rightarrow$}}
\end{picture}
\caption{Elements of $\Ch_8$ and their
multiplication.}\label{fig:f1}

\end{figure}

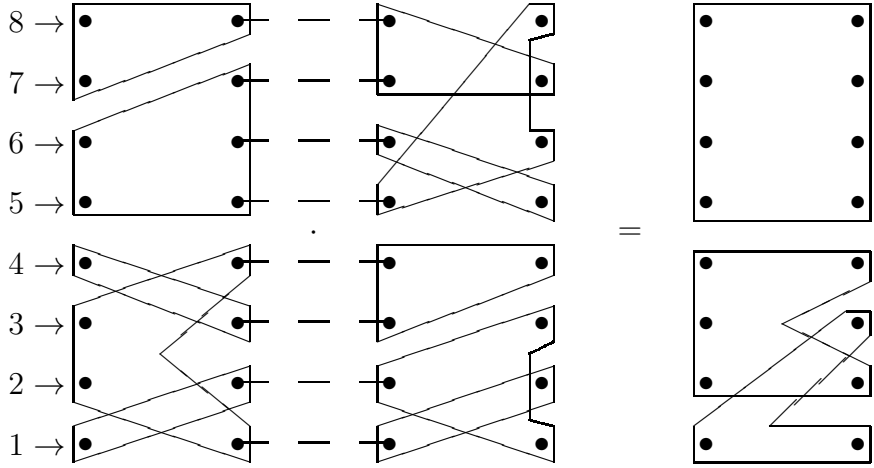
\begin{figure}
\special{em:linewidth 0.4pt} \unitlength 0.80mm
\linethickness{0.4pt}
\begin{picture}(150.00,75.00)
\put(20.00,00.00){\makebox(0,0)[cc]{$\bullet$}}
\put(20.00,10.00){\makebox(0,0)[cc]{$\bullet$}}
\put(20.00,20.00){\makebox(0,0)[cc]{$\bullet$}}
\put(20.00,30.00){\makebox(0,0)[cc]{$\bullet$}}
\put(20.00,40.00){\makebox(0,0)[cc]{$\bullet$}}
\put(20.00,50.00){\makebox(0,0)[cc]{$\bullet$}}
\put(20.00,60.00){\makebox(0,0)[cc]{$\bullet$}}
\put(20.00,70.00){\makebox(0,0)[cc]{$\bullet$}}
\put(45.00,10.00){\makebox(0,0)[cc]{$\bullet$}}
\put(45.00,20.00){\makebox(0,0)[cc]{$\bullet$}}
\put(45.00,30.00){\makebox(0,0)[cc]{$\bullet$}}
\put(45.00,40.00){\makebox(0,0)[cc]{$\bullet$}}
\put(45.00,50.00){\makebox(0,0)[cc]{$\bullet$}}
\put(45.00,60.00){\makebox(0,0)[cc]{$\bullet$}}
\put(45.00,70.00){\makebox(0,0)[cc]{$\bullet$}}
\put(45.00,00.00){\makebox(0,0)[cc]{$\bullet$}}
\drawline(18.00,-03.00)(18.00,03.00)
\drawline(18.00,03.00)(47.00,13.00)
\drawline(47.00,13.00)(47.00,07.00)
\drawline(47.00,07.00)(18.00,-03.00)
\drawline(18.00,07.00)(18.00,23.00)
\drawline(18.00,23.00)(47.00,33.00)
\drawline(47.00,33.00)(47.00,28.00)
\drawline(47.00,28.00)(32.25,15.00)
\drawline(32.25,15.00)(47.00,03.00)
\drawline(47.00,03.00)(47.00,-03.00)
\drawline(47.00,-03.00)(18.00,07.00)
\drawline(18.00,28.00)(18.00,33.00)
\drawline(18.00,33.00)(47.00,23.00)
\drawline(47.00,23.00)(47.00,17.00)
\drawline(47.00,17.00)(18.00,28.00)
\drawline(18.00,38.00)(18.00,52.00)
\drawline(18.00,52.00)(47.00,63.00)
\drawline(47.00,63.00)(47.00,38.00)
\drawline(47.00,38.00)(18.00,38.00)
\drawline(18.00,57.00)(18.00,73.00)
\drawline(18.00,73.00)(47.00,73.00)
\drawline(47.00,73.00)(47.00,68.00)
\drawline(47.00,68.00)(18.00,57.00)
\put(70.00,00.00){\makebox(0,0)[cc]{$\bullet$}}
\put(70.00,10.00){\makebox(0,0)[cc]{$\bullet$}}
\put(70.00,20.00){\makebox(0,0)[cc]{$\bullet$}}
\put(70.00,30.00){\makebox(0,0)[cc]{$\bullet$}}
\put(70.00,40.00){\makebox(0,0)[cc]{$\bullet$}}
\put(70.00,50.00){\makebox(0,0)[cc]{$\bullet$}}
\put(70.00,60.00){\makebox(0,0)[cc]{$\bullet$}}
\put(70.00,70.00){\makebox(0,0)[cc]{$\bullet$}}
\put(95.00,10.00){\makebox(0,0)[cc]{$\bullet$}}
\put(95.00,20.00){\makebox(0,0)[cc]{$\bullet$}}
\put(95.00,30.00){\makebox(0,0)[cc]{$\bullet$}}
\put(95.00,40.00){\makebox(0,0)[cc]{$\bullet$}}
\put(95.00,50.00){\makebox(0,0)[cc]{$\bullet$}}
\put(95.00,60.00){\makebox(0,0)[cc]{$\bullet$}}
\put(95.00,70.00){\makebox(0,0)[cc]{$\bullet$}}
\put(95.00,00.00){\makebox(0,0)[cc]{$\bullet$}}
\drawline(68.00,-03.00)(68.00,03.00)
\drawline(68.00,03.00)(97.00,13.00)
\drawline(97.00,13.00)(97.00,07.00)
\drawline(97.00,07.00)(68.00,-03.00)
\drawline(68.00,07.00)(68.00,13.00)
\drawline(68.00,13.00)(97.00,23.00)
\drawline(97.00,23.00)(97.00,17.00)
\drawline(97.00,17.00)(93.00,15.00)
\drawline(93.00,15.00)(93.00,04.00)
\drawline(93.00,04.00)(97.00,03.00)
\drawline(97.00,03.00)(97.00,-03.00)
\drawline(97.00,-03.00)(68.00,07.00)
\drawline(68.00,17.00)(68.00,33.00)
\drawline(68.00,33.00)(97.00,33.00)
\drawline(97.00,33.00)(97.00,28.00)
\drawline(97.00,28.00)(68.00,17.00)
\drawline(68.00,48.00)(68.00,53.00)
\drawline(68.00,53.00)(97.00,43.00)
\drawline(97.00,43.00)(97.00,37.00)
\drawline(97.00,37.00)(68.00,48.00)
\drawline(68.00,58.00)(68.00,73.00)
\drawline(68.00,73.00)(97.00,63.00)
\drawline(97.00,63.00)(97.00,58.00)
\drawline(97.00,58.00)(68.00,58.00)
\drawline(68.00,38.00)(68.00,43.00)
\drawline(68.00,43.00)(93.00,73.00)
\drawline(93.00,73.00)(97.00,73.00)
\drawline(97.00,73.00)(97.00,68.00)
\drawline(97.00,68.00)(93.00,67.00)
\drawline(93.00,67.00)(93.00,52.00)
\drawline(93.00,52.00)(97.00,52.00)
\drawline(97.00,52.00)(97.00,47.00)
\drawline(97.00,47.00)(68.00,38.00)
\drawline(45.00,00.30)(50.00,00.30)
\drawline(55.00,00.30)(60.00,00.30)
\drawline(65.00,00.30)(70.00,00.30)
\drawline(45.00,10.30)(50.00,10.30)
\drawline(55.00,10.30)(60.00,10.30)
\drawline(65.00,10.30)(70.00,10.30)
\drawline(45.00,20.30)(50.00,20.30)
\drawline(55.00,20.30)(60.00,20.30)
\drawline(65.00,20.30)(70.00,20.30)
\drawline(45.00,30.30)(50.00,30.30)
\drawline(55.00,30.30)(60.00,30.30)
\drawline(65.00,30.30)(70.00,30.30)
\drawline(45.00,40.30)(50.00,40.30)
\drawline(55.00,40.30)(60.00,40.30)
\drawline(65.00,40.30)(70.00,40.30)
\drawline(45.00,50.30)(50.00,50.30)
\drawline(55.00,50.30)(60.00,50.30)
\drawline(65.00,50.30)(70.00,50.30)
\drawline(45.00,60.30)(50.00,60.30)
\drawline(55.00,60.30)(60.00,60.30)
\drawline(65.00,60.30)(70.00,60.30)
\drawline(45.00,70.30)(50.00,70.30)
\drawline(55.00,70.30)(60.00,70.30)
\drawline(65.00,70.30)(70.00,70.30)
\drawline(120.00,37.00)(120.00,73.00)
\drawline(120.00,73.00)(149.00,73.00)
\drawline(149.00,73.00)(149.00,37.00)
\drawline(149.00,37.00)(120.00,37.00)
\drawline(120.00,08.00)(120.00,32.00)
\drawline(120.00,32.00)(149.00,32.00)
\drawline(149.00,32.00)(149.00,27.00)
\drawline(149.00,27.00)(134.50,20.00)
\drawline(134.50,20.00)(149.00,13.00)
\drawline(149.00,13.00)(149.00,08.00)
\drawline(149.00,08.00)(120.00,08.00)
\drawline(120.00,-03.00)(120.00,03.00)
\drawline(120.00,03.00)(145.00,22.00)
\drawline(145.00,22.00)(149.00,22.00)
\drawline(149.00,22.00)(149.00,18.00)
\drawline(149.00,18.00)(132.50,03.00)
\drawline(132.50,03.00)(149.00,03.00)
\drawline(149.00,03.00)(149.00,-03.00)
\drawline(149.00,-03.00)(120.00,-03.00)
\put(122.00,00.00){\makebox(0,0)[cc]{$\bullet$}}
\put(122.00,10.00){\makebox(0,0)[cc]{$\bullet$}}
\put(122.00,20.00){\makebox(0,0)[cc]{$\bullet$}}
\put(122.00,30.00){\makebox(0,0)[cc]{$\bullet$}}
\put(122.00,40.00){\makebox(0,0)[cc]{$\bullet$}}
\put(122.00,50.00){\makebox(0,0)[cc]{$\bullet$}}
\put(122.00,60.00){\makebox(0,0)[cc]{$\bullet$}}
\put(122.00,70.00){\makebox(0,0)[cc]{$\bullet$}}
\put(147.00,10.00){\makebox(0,0)[cc]{$\bullet$}}
\put(147.00,20.00){\makebox(0,0)[cc]{$\bullet$}}
\put(147.00,30.00){\makebox(0,0)[cc]{$\bullet$}}
\put(147.00,40.00){\makebox(0,0)[cc]{$\bullet$}}
\put(147.00,50.00){\makebox(0,0)[cc]{$\bullet$}}
\put(147.00,60.00){\makebox(0,0)[cc]{$\bullet$}}
\put(147.00,70.00){\makebox(0,0)[cc]{$\bullet$}}
\put(147.00,00.00){\makebox(0,0)[cc]{$\bullet$}}
\put(57.50,35.00){\makebox(0,0)[cc]{$\cdot$}}
\put(109.50,35.00){\makebox(0,0)[cc]{$=$}}

\put(12.00,00.00){\makebox(0,0)[cc]{$1\rightarrow$}}
\put(12.00,10.00){\makebox(0,0)[cc]{$2\rightarrow$}}
\put(12.00,20.00){\makebox(0,0)[cc]{$3\rightarrow$}}
\put(12.00,30.00){\makebox(0,0)[cc]{$4\rightarrow$}}
\put(12.00,40.00){\makebox(0,0)[cc]{$5\rightarrow$}}
\put(12.00,50.00){\makebox(0,0)[cc]{$6\rightarrow$}}
\put(12.00,60.00){\makebox(0,0)[cc]{$7\rightarrow$}}
\put(12.00,70.00){\makebox(0,0)[cc]{$8\rightarrow$}}
\end{picture}
\caption{Elements of $\IP_8$ and their
multiplication.}\label{fig:f2}

\end{figure}

Let $\IP_X$ be the subset of $\Ch_X$, containing those elements
$\bigl(A_i\bigr)_{i\in I}\in\Ch_X$ such that $A_i\cap
X\ne\varnothing$ and $A_i\cap X'\ne\varnothing$ for all $i\in I$.
Since the construction of $\Ch_X$, we have that $\IP_X$ is closed
under the multiplication in $\Ch_X$ and so $\IP_X$ is a subsemigroup
of $\Ch_X$. Observe that $\IP_X$ has the zero element, namely
$\{X\cup X'\}$. We will denote this element by $0$. Obviously, if
$\mid\!X\!\mid=\mid\!Y\!\mid$ then $\Ch_X\cong\Ch_Y$ and
$\IP_X\cong\IP_Y$. In the case when $X=\{1,\ldots,n\}$, it will be
convenient to denote $\Ch_X$ and $\IP_X$ by $\Ch_n$ and $\IP_n$
respectively. Figures~\ref{fig:f1} and~\ref{fig:f2} illustrate the
given notions for the case when $X=\{1,\ldots,8\}$, where we
consider elements of semigroups as couples of vertical rows of
points, divided into blocks. More precisely, the left vertical row
corresponds to the set $X$ and the right one to $X'$. The
multiplication $a\cdot b$ is just a gluing of elements $a$ and $b$
by dint of identifying the points of $X'$ from $a$ with the
corresponding elements of $X$ from $b$. On Fig.~\ref{fig:f1} we
present the equality
\begin{multline}
\bigl\{\{1,2,1'\},\{3,4\},\{5,2'\},\{3',4',5'\},\{6,7,6',7',8'\},\{8\}\bigr\}
\cdot\\
\bigl\{\{1,1'\},\{2,3,4\},\{2',3'\},\{5,5'\},\{6,4'\},\{7\},\{6',7'\},\{8,8'\}\bigr\}=\\
\bigl\{\{1,2,1'\},\{3,4\},\{2',3'\},\{5,5'\},\{6,7,4',8'\},\{6',7'\},\{8\}\bigr\}
\end{multline}
and on Fig.~\ref{fig:f2} we present the following one:
\begin{multline}
\bigl\{\{1,2'\},\{2,3,1',4'\},\{4,3'\},\{5,6,5',6',7'\},\{7,8,8'\}\bigr\}
\cdot\\
\bigl\{\{1,2'\},\{2,1',3'\},\{3,4,4'\},\{5,6',8'\},\{6,5'\},\{7,8,7'\}\bigr\}=\\
\bigl\{\{1,1',3'\},\{2,3,4,2',4'\},\{5,6,7,8,5',6',7',8'\}\bigr\}.
\end{multline}

Now we move to the proof of the fact that $\IP_X$ is isomorphic to
$\I^{\ast}_X$.

\section{$\IP_X$ is isomorphic to $\I^{\ast}_X$}\label{sec:Isomorphism}

The main goal of this section is to prove the following

\begin{theorem}\label{th:Isomorphism}
$\IP_X\cong\I^{\ast}_X$.
\end{theorem}

\begin{proof}
We begin with recalling one notion from~\cite{FL}. A \emph{block
bijection} of $X$ is a bijection between two quotient sets
$X/\sigma$ and $X/\tau$ for certain equivalence relations $\sigma$
and $\tau$ on $X$ such that $\lmod X/\sigma\rmod=\lmod X/\tau\rmod$.
We will need the following statement, stated in~\cite{FL} (one might
find it also in~\cite{Lawson}, Section 4.2).

\begin{lemma}[Lemma 2.1 from~\cite{FL}]\label{lm:FL}
If $\alpha$ is a biequivalence on $X$, then both
$\alpha\circ\alpha^{-1}$ and $\alpha^{-1}\circ\alpha$ are
equivalence relations on $X$. Moreover the map $\widetilde{\alpha}$
defined by $\widetilde{\alpha}:x(\alpha\circ\alpha^{-1})\mapsto
x\alpha$ for $x\in X$ is a block bijection of
$X/\alpha\circ\alpha^{-1}$ to $X/\alpha^{-1}\circ\alpha$.
Conversely, given equivalence relations $\beta$ and $\gamma$ on $X$
together with a block bijection $\mu:X/\beta\to X/\gamma$, a unique
biequivalence $\widehat{\mu}$ on $X$ inducing $\mu$ is given by:
$x\widehat{\mu}y$ if and only if $x\beta\mapsto y\gamma$ under the
block bijection $\mu$ (in which case
$\beta=\widehat{\mu}\circ\widehat{\mu}^{-1}$ and
$\gamma=\widehat{\mu}^{-1}\circ\widehat{\mu}$). Finally, the two
processes are reciprocal: $\widehat{\widetilde{\alpha}}=\alpha$ and
$\widetilde{\widehat{\mu}}=\mu$.
\end{lemma}
To define an isomorphism between $\IP_X$ and $\I^{\ast}_X$, we need
some auxiliary notation.

Let $a\in\IP_X$. Define the following relations $\rho_a$ and
$\lambda_a$ on $X$ as follows:
\begin{equation}
x\rho_{a}y~\mbox{if and only if}~x\equiv_a
y,~\mbox{and}~x\lambda_{a} y~\mbox{if and only if}~x'\equiv_a y',
\end{equation}
for $x,y\in X$. Since $\rho_{a}$ is a restriction of the relation
$\equiv_{a}$ to $X$, we obtain that $\rho_{a}$ is an equivalence
relation on $X$. From the definition of $\lambda_{a}$ and similar
arguments it follows that $\lambda_{a}$ is an equivalence relation
on $X$ as well. Remark that $a$ is not determined by $\lambda_{a}$
and $\rho_{a}$.

Define a map $\pi:\IP_X\to\I^{\ast}_X$ as follows: for $a\in\IP_X$
we put $\pi(a)=\widehat{\mu_a}$, where $\mu_a$ is a block bijection
from $X/\rho_a$ onto $X/\lambda_a$ such that the block $A$ of
$\rho_a$ is mapped under $\mu_a$ to that block $B$ of $\lambda_a$,
for which $A\cup B'$ is a block of $\equiv_a$. In view of our
definition of $\IP_X$ and Lemma~\ref{lm:FL}, we obtain that $\pi$ is
a bijection from $\IP_X$ onto $\I^{\ast}_X$.

We are left to prove that $\pi$ is a morphism from $\IP_X$ to
$\I^{\ast}_X$. Take $a,b\in\IP_X$. We need to prove that
$\widehat{\mu_{ab}}=\widehat{\mu_a}\widehat{\mu_b}=
\widehat{\mu_a}\circ\bigl(\widehat{\mu_a}^{-1}\circ\widehat{\mu_a}\vee
\widehat{\mu_b}\circ\widehat{\mu_b}^{-1}\bigr)\circ\widehat{\mu_b}$.
Notice that due to Lemma~\ref{lm:FL}, we have that
$\widehat{\mu_b}\circ\widehat{\mu_b}^{-1}=\rho_b$ and
$\widehat{\mu_a}^{-1}\circ\widehat{\mu_a}=\lambda_a$ and hence we
must establish that
$\widehat{\mu_{ab}}=\widehat{\mu_a}\circ\bigl(\lambda_a\vee
\rho_b\bigr)\circ\widehat{\mu_b}$. Note also that for all
$c\in\IP_X$ it follows immediately from the definition of $\mu_c$
that for all $x,y\in X$ one has $x\widehat{\mu_c}y$ if and only if
$x\equiv_c y'$. Finally, we recall that for equivalence relations
$\lambda$ and $\rho$ on $X$, the join $\lambda\vee\rho$ coincides
with the transitive closure of the relation $\lambda\cup\rho$.

Suppose firstly that $x\widehat{\mu_{ab}}y$, for some $x,y\in X$.
Then $x\equiv_{ab}y'$ and so there exist $c_1,\ldots$, $c_{2s-1}$,
$s\geq 1$, elements of $X$, such that $x\equiv_{a} c'_1$,
$c_1\equiv_{b} c_2$, $c'_2\equiv_{a} c'_3,\ldots,$
$c'_{2s-2}\equiv_{a} c'_{2s-1}$, and $c_{2s-1}\equiv_{b} y'$. Then
we have $x\widehat{\mu_{a}}c_1$, $c_1\rho_{b} c_2$, $c_2\lambda_{a}
c_3,\ldots,$ $c_{2s-2}\lambda_{a} c_{2s-1}$, and
$c_{2s-1}\widehat{\mu_{b}}y$. Thus, we have $x\widehat{\mu_{a}}c_1$,
$c_1(\lambda_a\vee\rho_b)c_{2s-1}$ and $c_{2s-1}\widehat{\mu_{b}}y$,
whence $(x,y)\in\widehat{\mu_a}\circ\bigl(\lambda_a\vee
\rho_b\bigr)\circ\widehat{\mu_b}$.

Conversely, suppose that
$(x,y)\in\widehat{\mu_a}\circ\bigl(\lambda_a\vee
\rho_b\bigr)\circ\widehat{\mu_b}$. Then there exist $c,d\in X$ such
that $x\widehat{\mu_{a}}c$, $c(\lambda_a\vee\rho_b)d$ and
$d\widehat{\mu_{b}}y$. Then we have $x\equiv_a c'$ and $d\equiv_b
y'$. Notice that if $c\lambda_a r$ then $x\equiv_a r'$ and if
$t\rho_b d$ then $t\equiv_b y'$. Hence, taking to account
$c(\lambda_a\vee\rho_b)d$, there exist $c_1,\ldots$, $c_{2s-1}$,
$s\geq 1$, elements of $X$, such that $x\equiv_{a}c_1'$,
$c_1\rho_{b} c_2$, $c_2\lambda_{a} c_3,\ldots,$ $c_{2s-2}\lambda_{a}
c_{2s-1}$, and $c_{2s-1}\equiv_{b}y'$. These imply $x\equiv_{a}
c'_1$, $c_1\equiv_{b} c_2$, $c'_2\equiv_{a} c'_3,\ldots,$
$c'_{2s-2}\equiv_{a} c'_{2s-1}$, and $c_{2s-1}\equiv_{b} y'$. Thus
$x\equiv_{ab}y'$, whence $x\widehat{\mu_{ab}}y$.

The proof of the theorem is complete.
\end{proof}

As a consequence of Theorem~\ref{th:Isomorphism} we obtain the
following statement.

\begin{proposition} \label{cor:IP-inverse}
$\IP_X$ is an inverse semigroup.
\end{proposition}

\begin{proof}
Follows from the fact that $\I^{\ast}_X$ is inverse, see~\cite{FL}.
\end{proof}

Due to what we have already obtained, we can now call $\IP_X$ the
\emph{inverse partition semigroup} on the set $X$.

\section{Green's relations and the natural order in $\IP_X$}\label{sec:green}

We begin this section with description of Green's relations on
$\IP_X$. But before we need some preparation.

First notice that it follows immediately from the definition of
multiplication in $\IP_X$ that
\begin{equation}\label{eq:remark-for-ranks}
\rho_{ab}\supseteq\rho_a~\mbox{and}~\lambda_{ab}\supseteq\lambda_b~\mbox{for
all}~a,b\in\IP_X.
\end{equation}
Then we obtain that every $\rho_{ab}$-class is a union of some
$\rho_a$-classes and that every $\lambda_{ab}$-class is a union of
some $\lambda_b$-classes.

Note also that the cardinal number of the set of all
$\rho_a$-classes and the cardinal number of the set of all
$\lambda_a$-classes coincide with the cardinal number of the set of
all $\equiv_a$-classes. Denote this common number by $\rank(a)$. We
will call this number the \emph{rank} of $a$. Due
to~\eqref{eq:remark-for-ranks}, we have
\begin{equation} \label{eq:ranks}
\rank(ab)\leq\min\bigl\{\rank(a),\rank(b)\bigr\}~\mbox{for
all}~a,b\in\IP_X.
\end{equation}
Note that if $a=\bigl(A_i\cup B_i'\bigr)_{i\in I}$ then
$\rank(a)=\lmod I\rmod$. We denote the Green's relations in the
standard way: $\GR$, $\GL$, $\GH$, $\GD$, and $\GJ$
(see~\cite{Howie}).

\begin{theorem} \label{th:Green}
Let $a,b\in\IP_X$. Then
\begin{enumerate}
\item\label{R}
$a\GR b$ if and only if $\rho_a=\rho_b$;
\item\label{L}
$a\GL b$ if and only if $\lambda_a=\lambda_b$;
\item\label{H}
$a\GH b$ if and only if $\rho_a=\rho_b$ and $\lambda_a=\lambda_b$
hold simultaneously;
\item\label{J=D}
$a\GJ b$ if and only if $a\GD b$ if and only if $\rank(a)=\rank(b)$;
\item\label{card-sem}
$\mid\!\IP_n\!\mid=\sum\limits_{k=1}^{n}\bigl(s(n,k)\bigr)^2\cdot
k!$, where $s(n,k)$ denotes the Stirling number of the second kind;
\item\label{card-idem}
$\mid\!E(\IP_n)\!\mid=\mathrm{B}_n$, where $\mathrm{B}_n$ denotes
the Bell number.
\end{enumerate}
\end{theorem}

\begin{proof}
In view of Theorem~\ref{th:Isomorphism}, these statements are just
reformulations of those of Theorem 2.2 from~\cite{FL}.
\end{proof}

Now we move to description of the group of units of $\IP_X$. Denote
by $\mathcal{S}_X$ the \emph{symmetric group} on $X$. Set a map
$\eta:\mathcal{S}_X\to\IP_X$ as follows:
\begin{equation}
\eta(g)=\bigl(\{x,g(x)'\}\bigr)_{x\in X}~\mbox{for
all}~g\in\mathcal{S}_X.
\end{equation}

\begin{lemma} \label{lm:S_XembedsintoIP_X}
The map $\eta$ is an injective homomorphism.
\end{lemma}

\begin{proof}
That $\eta$ is a homomorphism, follows from the definition of the
multiplication in $\IP_X$. If now $\eta(g_1)=\eta(g_2)$ for some
$g_1,g_2\in\mathcal{S}_X$, then $g_1(x)=g_2(x)$ for all $x\in X$ and
so $g_1=g_2$. This completes the proof.
\end{proof}

As a consequence of Lemma~\ref{lm:S_XembedsintoIP_X} we obtain that
$\IP_X$ contains a subgroup $\eta(\mathcal{S}_X)$, isomorphic to
$\mathcal{S}_X$. Let us identify this subgroup with $\mathcal{S}_X$.
Clearly, the identity element $1$ of $\Sym_X$ is the identity
element of $\IP_X$. Using Theorem~\ref{th:Green}, we obtain now the
following corollary.

\begin{proposition} \label{cor:Invertible-Elements}
The group of all invertible elements of $\IP_X$ coincides with
$\Sym_X$.
\end{proposition}

\begin{proof}
Since the maximal subgroup of an arbitrary semigroup coincides with
some $\GH$-class of this semigroup (see~\cite{Howie}), we obtain
that an element $g$ is invertible in $\IP_X$ if and only if $g\GH
1$. Due to Theorem~\ref{th:Green}, this is equivalent to
$g\in\Sym_X$.
\end{proof}

Let us now switch to the description of the natural order on
$\IP_X$. But before, we need to describe the idempotents of $\IP_X$.

\begin{lemma} \label{pr:idempotentIP}
Let $e\in\IP_X$. Then $e$ is an idempotent if and only if there is a
partition $X=\bigcup\limits_{i\in I}^{\cdot}E_i$ such that
$e=\bigl(E_i\cup E_i'\bigr)_{i\in I}$. In addition, for idempotents
$e$ and $f$ the elements $ef$ and $fe$ coincide with the minimum
equivalence relation on $X\cup X'$, which contains $e$ and $f$.
\end{lemma}

\begin{proof}
Let us prove firstly the first part of the statement. The
sufficiency of it is obvious.

Let now $e$ be an idempotent of $\IP_X$. Let $A\cup B'$ be some
block in $e$. Suppose that $A\setminus B\ne \varnothing$. Then there
is $a\in A$ such that $a\notin B$. Take an arbitrary $b$ of $B$.
Take also $c\in X$ such that $c\equiv_{e}a'$. Then $c\notin A$.
Indeed, otherwise we would have $a\equiv_{e}c\equiv_{e}a'$ which
implies $a\in B$. Thus, $c\notin A$.

Now due to $c\equiv_{e}a'$ and $a\equiv_{e}b'$, we obtain that
$c\equiv_{e^2}b'$. But the latter gives us $c\in A$. We get a
contradiction. Thus, $A\setminus B=\varnothing$ and so $A\subseteq
B$. Analogously, $B\subseteq A$. Thus, every block of $e$ has the
form $A\cup A'$ for certain $A\subseteq X$. This completes the proof
of the first part of the statement. The second one now follows
immediately from the definition of the multiplication in $\IP_X$.
\end{proof}

\begin{proposition}\label{pr:Omega}
Let $a,b\in\IP_X$. Then $a\leq b$ if and only if
$\equiv_a\supseteq\equiv_b$.
\end{proposition}

\begin{proof}
Let $a=\bigl(A_i\cup B_i'\bigr)_{i\in I}$ and $b=\bigl(C_j\cup
D_j'\bigr)_{j\in J}$.

Suppose first that $\equiv_{b}\subseteq\equiv_{a}$. Then we have
that for all $i\in I$, $A_i\cup B_i'$ is a union of some blocks
$C_j\cup D_j'$, $j\in J$. Put $f=\bigl(B_i\cup B_i'\bigr)_{i\in I}$.
Then we obtain that $a=bf$. It remains to note that, due to
Lemma~\ref{pr:idempotentIP}, $f$ is an idempotent.

Suppose now that there is an idempotent $e$ of $\IP_X$ such that
$a=be$. Due to Lemma~\ref{pr:idempotentIP}, we have that
$e=\bigl(E_k\cup E_k'\bigr)_{k\in K}$ for some partition
$X=\bigcup\limits_{k\in K}^{\cdot}E_k$. Take now $(x,y)\in\equiv_b$.
There is $z$ of $X$ such that $z'$ is $\equiv_b$--equivalent to $x$
and $y$. Then, since $z\equiv_e z'$, we obtain that
$(x,y)\in\equiv_{be}$ or just that $(x,y)\in\equiv_{a}$. This
completes the proof.
\end{proof}

Now we are able to characterize the trace of $\IP_X$.

\begin{proposition}\label{pr:trace}
Let $a,b\in\tr(\IP_X)$. The product $a\ast b$ is defined if
$\lambda_{a}=\rho_{b}$ and in this case
$\pi(a)\circ\pi(b)\in\I^{\ast}_X$ and $a\ast
b=\pi^{-1}(\pi(a)\circ\pi(b))$.
\end{proposition}

\begin{proof}
It is known that for $x,y\in\tr(S)$, where $S$ is an inverse
semigroup, the product $x\ast y$ is defined if and only if
$x^{-1}x=yy^{-1}$ (see~\cite{Meakin}). Note also that, using
Lemma~\ref{pr:idempotentIP}, we have that for every $x\in\IP_X$ the
condition $\rho_{x}=\lambda_{x}$ holds if and only if $x\in
E(\IP_X)$. In addition, for $e,f\in E(\IP_X)$ we have that
$\lambda_{e}=\lambda_{f}$ if and only if $e=f$. Hence, $a\ast b$ is
defined if and only if $a^{-1}a=bb^{-1}$ if and only if
$\lambda_{a^{-1}a}=\rho_{bb^{-1}}$. It remains to notice that since
$a^{-1}a\GL a$ and $bb^{-1}\GR b$, using Theorem~\ref{th:Green}, we
have $\lambda_{a^{-1}a}=\lambda_{a}$ and $\rho_{bb^{-1}}=\rho_{b}$.

If now $a\ast b$ is defined then $\pi(a)\ast\pi(b)$ is defined in
$\I^{\ast}_X$ and then $\pi(a)\ast\pi(b)=\pi(a)\circ\pi(b)$
(see~\cite{Lawson}). The statement follows.
\end{proof}

The following proposition is concerned with $\im(\IP_X)$, the
imprint of $\IP_X$.

\begin{proposition}
Let $e\in E(\IP_X)$ and $a\in\IP_X$. The product $e\star a$ is
defined if and only if $\rho_{a}\subseteq\rho_{e}$.
\end{proposition}

\begin{proof}
By the definition of imprint, we have that $e\star a$ is defined if
and only if $e\leq aa^{-1}$, which, in view of
Proposition~\ref{pr:Omega}, holds if and only if
$\equiv_{aa^{-1}}\subseteq\equiv_{e}$ which is equivalent to
$\rho_{aa^{-1}}\subseteq\rho_{e}$. It remains to notice that
$\rho_{a}=\rho_{aa^{-1}}$.
\end{proof}

\section{Generating set, ideals and maximal subsemigroups of $\IP_n$} \label{sec:some-objects}

To begin this section, we put some auxiliary notations. Let
$A\subseteq X$. Define an element $\tau_{A}$ of $\IP_X$ as follows:
\begin{equation}
\tau_A=\bigl\{A\cup A',\{x,x'\}_{x\in X\setminus A}\bigr\}.
\end{equation}
Clearly, $\tau_X$ is the zero element of $\IP_X$. If $x$ and $y$ are
distinct elements of $X$, we will use the notation
$\tau_{x,y}=\tau_{\{x,y\}}$.

Suppose that $\mid\!X\!\mid\geq 3$. For pairwise distinct elements
$x,y,z$ of $X$ define an element $\xi_{x,y,z}$ as follows:
\begin{equation}
\xi_{x,y,z}=\bigl\{\{x,y,x'\},\{z,y',z'\},\{t,t'\}_{t\in X\setminus
\{x,y,z\}}\bigr\}.
\end{equation}
If necessary, we will write $\xi_{x,y,z}^{X}$ instead of
$\xi_{x,y,z}$ to stress on that $\xi_{x,y,z}\in\IP_X$.

\begin{lemma}\label{lm:XI-and-TAU}
Let $\mid\!X\!\mid\geq 3$. Then
\begin{multline}
g^{-1}\xi_{x,y,z}g=
\xi_{g(x),g(y),g(z)},~g^{-1}\tau_{x,y}g=\tau_{g(x),g(y)},\\
\xi_{x,y,z}^{2}=\tau_{\{x,y,z\}}~\mbox{and}~\xi_{x,y,z}\xi_{z,y,x}=\tau_{x,y}
\end{multline}
for all pairwise distinct $x,y,z\in X$ and $g\in\Sym_X$.
\end{lemma}

\begin{proof}
Direct calculation.
\end{proof}

Now our local goal is to provide a generating set for $\IP_n$ (see
Proposition~\ref{pr:Generating-SystemsIP}). In order to do this we
will construct an inverse subsemigroup $\IT_n$ of $\IP_n$ (see
below), which is interesting itself as a semigroup. In addition, the
notion of $\IT_n$ will help us to describe all the maximal
subsemigroups of $\IP_n$. So we are starting with putting some
auxiliary notations.

Let $n\geq 2$. Set $\IT_n=\langle \Sym_n,\tau_{1,2}\rangle$. Set
also $\IT_1=\IP_1$. Let $\rho$ be some equivalence relation on
$\{1,\ldots,n\}$. Define a \emph{type} of the relation $\rho$ as a
tuple $(t_1,\ldots,t_n)$, where $t_i$ denotes the number of all
$i$-element $\rho$-classes, $1\leq i\leq n$. The following
proposition shows that $\IT_n$ is an inverse subsemigroup of
$\IP_n$. But before, we give one more definition: an element $a$ of
$\IP_n$ is said to be \emph{special} if
\begin{equation}
x\equiv_a
y'~\mbox{implies}~\mid\!x\rho_a\!\mid=\mid\!y\lambda_a\!\mid~\mbox{for
all}~x,y\in\{1,\ldots,n\}.
\end{equation}

\begin{proposition}\label{pr:IT_n}
The following statements hold:
\begin{enumerate}
\item\label{1}
$\IT_n$ is an inverse subsemigroup of $\IP_n$;
\item\label{2}
$\tau_A\in \IT_n$ for all $A\subseteq\{1,\ldots,n\}$;
\item\label{3}
the elements of $\IT_n$ are precisely all special elements of
$\IP_n$;
\item\label{4}
if $a\in\IT_n$ then the types of $\rho_a$ and $\lambda_a$ coincide.
\end{enumerate}
\end{proposition}

\begin{proof}
We will assume that $n\geq 2$ as all the statements hold in the case
when $n=1$.

Since $\Sym_n$ is a subgroup of $\IP_n$ and $\tau_{1,2}$ is an
idempotent in $\IP_n$, we obtain that $\IT_n$ is an inverse
subsemigroup of $\IP_n$. This completes the proof of~\ref{1}).

Note that, due to Lemma~\ref{lm:XI-and-TAU}, we have that
$\tau_{x,y}\in\IT_n$ for all distinct $x$ and $y$ of
$\{1,\ldots,n\}$. Now the statement~\ref{2}) follows from the
equality $\tau_{\{x\}}=1$, for all $x\in\{1,\ldots,n\}$, and the
fact that if $A=\{x_1,\ldots,x_k\}$, $k\geq 2$, then
\begin{equation}
\tau_{A}=\prod\limits_{i=1}^{k-1}\tau_{x_k,x_{k+1}}.
\end{equation}

Let us prove~\ref{3}). Let $a=\bigl(A_i\cup B_i'\bigr)_{i\in I}$ be
an element of $\IP_n$ such that $x\equiv_a y'$ implies
$\mid\!x\rho_a\!\mid=\mid\!y\lambda_a\!\mid$ for all
$x,y\in\{1,\ldots,n\}$. Then $\mid\!A_i\!\mid=\mid\!B_i\!\mid$ for
all $i\in I$ and so there exists $g\in\Sym_n$ such that
$ga=\bigl(B_i\cup B_i'\bigr)_{i\in I}$. Now due to~\ref{2}), we have
that
\begin{equation}
a=g^{-1}\cdot\prod\limits_{i\in I}\tau_{B_i}\in\IT_n.
\end{equation}
Conversely, suppose that $a\in\IT_n$. Note that $\tau_{1,2}$ is
special and all the elements of $\Sym_n$ are special, too. Hence, to
prove that $a$ is special, it is enough to prove that if $b\in\IP_n$
is special then $b\tau_{1,2}$ is special and $bg$ is special for all
$g\in\Sym_n$. Suppose that $b=\bigl(C_i\cup D_i'\bigr)_{i\in
K}\in\IP_n$ is special. Then, obviously, $bg$ is also special for
all $g\in\Sym_n$. We have two cases.

\emph{Case 1}. There is $i\in K$ such that $D_i\supseteq\{1,2\}$.
Then $b\tau_{1,2}=b$ is special.

\emph{Case 2}. There are distinct $i$ and $j$ of $K$ such that $1\in
D_i$ and $2\in D_j$. Then $b\tau_{1,2}=\bigl\{\bigl(C_i\cup
C_j\bigr)\bigcup\bigl(D_i\cup D_j\bigr)',\bigl(C_k\cup
D_k'\bigr)_{k\in K\setminus\{i,j\}}\bigr\}$ is, obviously, special.
This completes the proof of~\ref{3}).

The statement~\ref{4}) follows immediately from~\ref{3}).
\end{proof}

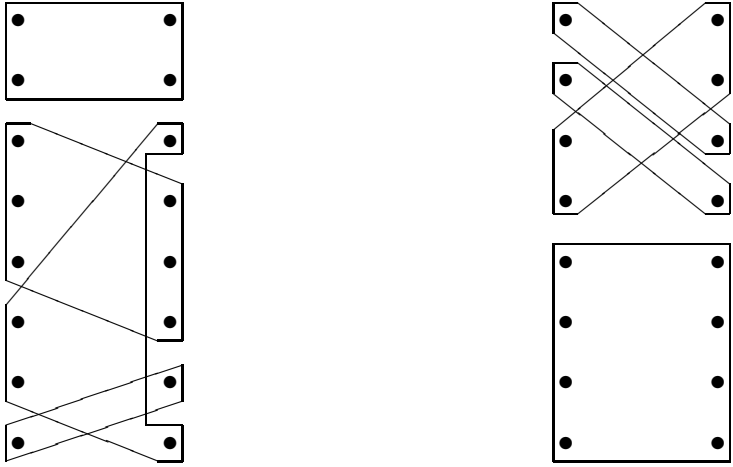
\begin{figure}
\special{em:linewidth 0.4pt} \unitlength 0.80mm
\linethickness{0.4pt}
\begin{picture}(150.00,75.00)
\put(15.00,10.00){\makebox(0,0)[cc]{$\bullet$}}
\put(15.00,20.00){\makebox(0,0)[cc]{$\bullet$}}
\put(15.00,30.00){\makebox(0,0)[cc]{$\bullet$}}
\put(15.00,40.00){\makebox(0,0)[cc]{$\bullet$}}
\put(15.00,50.00){\makebox(0,0)[cc]{$\bullet$}}
\put(15.00,60.00){\makebox(0,0)[cc]{$\bullet$}}
\put(15.00,70.00){\makebox(0,0)[cc]{$\bullet$}}
\put(15.00,00.00){\makebox(0,0)[cc]{$\bullet$}}

\put(40.00,00.00){\makebox(0,0)[cc]{$\bullet$}}
\put(40.00,10.00){\makebox(0,0)[cc]{$\bullet$}}
\put(40.00,20.00){\makebox(0,0)[cc]{$\bullet$}}
\put(40.00,30.00){\makebox(0,0)[cc]{$\bullet$}}
\put(40.00,40.00){\makebox(0,0)[cc]{$\bullet$}}
\put(40.00,50.00){\makebox(0,0)[cc]{$\bullet$}}
\put(40.00,60.00){\makebox(0,0)[cc]{$\bullet$}}
\put(40.00,70.00){\makebox(0,0)[cc]{$\bullet$}}

\put(105.00,10.00){\makebox(0,0)[cc]{$\bullet$}}
\put(105.00,20.00){\makebox(0,0)[cc]{$\bullet$}}
\put(105.00,30.00){\makebox(0,0)[cc]{$\bullet$}}
\put(105.00,40.00){\makebox(0,0)[cc]{$\bullet$}}
\put(105.00,50.00){\makebox(0,0)[cc]{$\bullet$}}
\put(105.00,60.00){\makebox(0,0)[cc]{$\bullet$}}
\put(105.00,70.00){\makebox(0,0)[cc]{$\bullet$}}
\put(105.00,00.00){\makebox(0,0)[cc]{$\bullet$}}

\put(130.00,00.00){\makebox(0,0)[cc]{$\bullet$}}
\put(130.00,10.00){\makebox(0,0)[cc]{$\bullet$}}
\put(130.00,20.00){\makebox(0,0)[cc]{$\bullet$}}
\put(130.00,30.00){\makebox(0,0)[cc]{$\bullet$}}
\put(130.00,40.00){\makebox(0,0)[cc]{$\bullet$}}
\put(130.00,50.00){\makebox(0,0)[cc]{$\bullet$}}
\put(130.00,60.00){\makebox(0,0)[cc]{$\bullet$}}
\put(130.00,70.00){\makebox(0,0)[cc]{$\bullet$}}

\drawline(13.00,-03.00)(13.00,03.00)
\drawline(13.00,03.00)(42.00,13.00)
\drawline(42.00,13.00)(42.00,07.00)
\drawline(42.00,07.00)(13.00,-03.00)

\drawline(13.00,57.00)(13.00,73.00)
\drawline(13.00,73.00)(42.00,73.00)
\drawline(42.00,73.00)(42.00,57.00)
\drawline(42.00,57.00)(13.00,57.00)

\drawline(13.00,27.00)(13.00,53.00)
\drawline(13.00,53.00)(17.00,53.00)
\drawline(17.00,53.00)(42.00,43.00)
\drawline(42.00,43.00)(42.00,17.00)
\drawline(42.00,17.00)(38.00,17.00)
\drawline(38.00,17.00)(13.00,27.00)

\drawline(13.00,07.00)(13.00,23.00)
\drawline(13.00,23.00)(38.00,53.00)
\drawline(38.00,53.00)(42.00,53.00)
\drawline(42.00,53.00)(42.00,48.00)
\drawline(42.00,48.00)(36.00,48.00)
\drawline(36.00,48.00)(36.00,03.00)
\drawline(36.00,03.00)(42.00,03.00)
\drawline(42.00,03.00)(42.00,-03.00)
\drawline(42.00,-03.00)(38.00,-03.00)
\drawline(38.00,-03.00)(13.00,07.00)

\drawline(103.00,-03.00)(103.00,33.00)
\drawline(103.00,33.00)(132.00,33.00)
\drawline(132.00,33.00)(132.00,-03.00)
\drawline(132.00,-03.00)(103.00,-03.00)

\drawline(103.00,58.00)(103.00,63.00)
\drawline(103.00,63.00)(107.00,63.00)
\drawline(107.00,63.00)(132.00,43.00)
\drawline(132.00,43.00)(132.00,38.00)
\drawline(132.00,38.00)(128.00,38.00)
\drawline(128.00,38.00)(103.00,58.00)

\drawline(103.00,68.00)(103.00,73.00)
\drawline(103.00,73.00)(107.00,73.00)
\drawline(107.00,73.00)(132.00,53.00)
\drawline(132.00,53.00)(132.00,48.00)
\drawline(132.00,48.00)(128.00,48.00)
\drawline(128.00,48.00)(103.00,68.00)

\drawline(103.00,38.00)(103.00,52.00)
\drawline(103.00,52.00)(128.00,73.00)
\drawline(128.00,73.00)(132.00,73.00)
\drawline(132.00,73.00)(132.00,58.00)
\drawline(132.00,58.00)(107.00,38.00)
\drawline(107.00,38.00)(103.00,38.00)

\end{picture}
\caption{Elements of $\IT_8$.}\label{fig:f3}

\end{figure}

As a consequence of~\ref{4}) of Proposition~\ref{pr:IT_n}, we can
now call $\IT_n$ the \emph{inverse type-preserving semigroup} of
degree $n$. We give an illustration of elements of $\IT_8$ on
Fig.~\ref{fig:f3}. It also follows from Proposition~\ref{pr:IT_n}
that $\IT_n=\Sym_nE(\IP_n)$, that is $\IT_n$ is the greatest
factorizable inverse submonoid of $\IP_n$. Remark that $\IT_n$ (more
precisely, $\pi(\IT_n)$, the greatest factorizable inverse submonoid
of $\I^{\ast}_X$) appeared in~\cite{F},~\cite{FL} and~\cite{AO}
under the name of the \emph{monoid of uniform block permutations}.

The following proposition gives us an example of a generating system
of $\IP_n$. But to prove this proposition, we need some auxiliary
facts.

\begin{lemma}\label{lm:BeforePrGeneratingSystems}
Let $n\geq 3$, $a\in\IP_n$ and $\rank(a)=n-1$. Then either
$a\in\xi_{x,y,z}\Sym_n$ or $a\in\tau_{x,y}\Sym_n$ for some pairwise
distinct $x,y,z\in\{1,\ldots,n\}$.
\end{lemma}

\begin{proof}
Straightforward.
\end{proof}

Take $n\in\mathbb{N}$. Set
$\Pi_n=\bigl\{q\in\IP_{n+1}:~q~\mbox{contains the
block}~\{n+1,(n+1)'\}\bigr\}$.

\begin{lemma}\label{lm:Pi_n}
Let $n\in\mathbb{N}$. Then the map $a\mapsto a\cup\{n+1,(n+1)'\}$,
$a\in\IP_n$, is an isomorphism from $\IP_n$ onto $\Pi_n$, which maps
$\xi_{1,2,3}^{\{1,\ldots,n\}}$ to $\xi_{1,2,3}^{\{1,\ldots,n+1\}}$.
\end{lemma}

\begin{proof}
Obvious.
\end{proof}

\begin{proposition}\label{pr:Generating-SystemsIP}
Let $n\geq 3$. Then $\IP_n=\langle \Sym_n, \xi_{1,2,3}\rangle$.
Moreover, for $u\in\IP_n$, $\IP_n=\langle \Sym_n,u\rangle$ if and
only if $u\in \Sym_n\xi_{1,2,3}\Sym_n$.
\end{proposition}

\begin{proof}
We will prove the statement that $\IP_n=\langle \Sym_n,
\xi_{1,2,3}\rangle$ for all $n\geq 3$ by the complete induction on
$n$.

First, let us verify that the basis of the induction, the case when
$n=3$, holds. We are to prove that $\IP_3=\langle \Sym_3,
\xi_{1,2,3}\rangle$. Note that, due to Lemma~\ref{lm:XI-and-TAU},
$0=\xi_{1,2,3}^2$. Thus, we are left to prove that every element $v$
of $\IP_3$ such that $\rank(v)=2$, belongs to $\langle \Sym_3,
\xi_{1,2,3}\rangle$. But this follows from Lemmas
\ref{lm:XI-and-TAU} and~\ref{lm:BeforePrGeneratingSystems}. Thus,
the basis of induction holds.

Assume now that the proposition of induction holds for all numbers
$k$, $3\leq k\leq n$. We are going to prove now that
$\IP_{n+1}=\langle \Sym_{n+1}, \xi_{1,2,3}\rangle$. Let $a\in
\IP_{n+1}$. Then there is $g\in \Sym_{n+1}$ such that $b=ag$
contains a block $\bigl(E\cup \{n+1\}\bigr)\bigcup \bigl(F\cup
{\{n+1\}}\bigr)'$ for certain subsets $E$ and $F$ of
$\{1,\ldots,n\}$. Note that, due to Lemma~\ref{lm:XI-and-TAU},
$\tau_{x,y}$ and $\xi_{x,y,z}$ are both elements of $\langle
\Sym_{n+1}, \xi_{1,2,3}\rangle$ for all pairwise distinct
$x,y,z\in\{1,\ldots,n\}$. Then taking to account Proposition
\ref{pr:IT_n}, we obtain that $\IT_n\subseteq\langle \Sym_{n+1},
\xi_{1,2,3}\rangle$. In particular, $0\in\langle \Sym_{n+1},
\xi_{1,2,3}\rangle$. Thus, without loss of generality we may assume
that $a\ne 0$, which implies $b\ne 0$. Suppose that all the blocks
of $b$, except $\bigl(E\cup \{n+1\}\bigr)\bigcup \bigl(F\cup
{\{n+1\}}\bigr)'$, are precisely $E_i\cup F_i'$, $1\leq i\leq k$. By
the proposition of induction and Lemma~\ref{lm:Pi_n}, we obtain that
\begin{equation}
c=\bigl\{\bigl(E\cup E_1\bigr)\cup \bigl(F\cup F_1\bigr)',E_2\cup
{F_2}',\ldots,E_k\cup {F_k}',\{n+1,(n+1)'\}\bigr\}
\end{equation}
is an element of $\langle \Sym_{n+1}, \xi_{1,2,3}\rangle$. We have
four possibilities.

\emph{Case 1}. $E=\varnothing$ and $F=\varnothing$. Then
$b=c\in\langle \Sym_{n+1}, \xi_{1,2,3}\rangle$.

\emph{Case 2}. $E=\varnothing$ and
$F=\{f_1,\ldots,f_m\}\ne\varnothing$. Fix an element $f\in F_1$.
Then $b=c\cdot\prod\limits_{i=1}^{m}\xi_{f,f_i,n+1}$ and so
$b\in\langle \Sym_{n+1}, \xi_{1,2,3}\rangle$.

\emph{Case 3}. $E=\{e_1,\ldots,e_l\}\ne\varnothing$ and
$F=\varnothing$. Fix an element $e\in E_1$. Then
$b=\prod\limits_{i=1}^{l}\xi_{n+1,f_i,e}\cdot c$, whence
$b\in\langle \Sym_{n+1}, \xi_{1,2,3}\rangle$.

\emph{Case 4}. $E\ne\varnothing$ and $F\ne\varnothing$. Put
$d=\bigl\{E\cup F',E_1\cup {F_1}',\ldots,E_k\cup
{F_k}',\{n+1,(n+1)'\}\bigr\}$. Due to proposition of induction and
Lemma~\ref{lm:Pi_n}, we have that $d\in\langle \Sym_{n+1},
\xi_{1,2,3}\rangle$. Then $b=\tau_{E\cup\{n+1\}}
d\tau_{F\cup\{n+1\}}\in\langle \Sym_{n+1}, \xi_{1,2,3}\rangle$.

In all these cases we obtained that $b$ belongs to $\langle
\Sym_{n+1}, \xi_{1,2,3}\rangle$ and so does $a$.

Thus, we have just proved that $\IP_n=\langle\Sym_{n},
\xi_{1,2,3}\rangle$ for all $n\geq 3$. This implies that
$\IP_n=\langle \Sym_{n},u\rangle$ for all
$u\in\Sym_{n}\xi_{1,2,3}\Sym_{n}$. Conversely, suppose that
$\IP_n=\langle \Sym_{n},u\rangle$ for some $u\in\IP_n$. Then, due to
\eqref{eq:ranks}, we obtain that $\rank(u)=n-1$. Now taking to
account Lemmas~\ref{lm:BeforePrGeneratingSystems} and
\ref{lm:XI-and-TAU}, we have that either
$u\in\Sym_n\xi_{1,2,3}\Sym_n$ or $u\in\Sym_n\tau_{1,2}\Sym_{n}$. But
$u\in\Sym_n\tau_{1,2}\Sym_{n}$ is impossible. Indeed, otherwise we
would have $\langle\Sym_n,\xi_{1,2,3}\rangle=\IT_n$ and it remains
to note that, due to~\ref{3}) of Proposition~\ref{pr:IT_n},
$\xi_{1,2,3}\notin\IT_n$ when $n\geq 3$. Hence,
$u\in\Sym_n\xi_{1,2,3}\Sym_n$ holds, as was required. This completes
the proof.
\end{proof}

Let $k\in\mathbb{N}$, $k\leq n$. Set
$I_k=\bigl\{a\in\IP_n:~\rank(a)\leq k\bigr\}$. Note that
\begin{equation}\label{eq:Ideal-Chain}
\{0\}=I_1\subset I_2\subset\ldots\subset I_n=\IP_n.
\end{equation}
We will prove in the following proposition that these sets exhaust
all the double-sided ideals (or just ideals) of $\IP_n$.

\begin{proposition}\label{pr:Ideals}
Let $I$ be an ideal of $\IP_n$ and $k\in\mathbb{N}$ such that $k\leq
n$. Then
\begin{enumerate}
\item\label{Ideals-1}
for all $b\in\IP_n$, $I_k=\IP_nb\IP_n$ if and only if $\rank(b)=k$;
\item\label{Ideals-2}
$I=I_m$ for some $m\in\mathbb{N}$, $m\leq n$;
\item\label{Ideals-3}
$I=\IP_na\IP_n$ for certain $a\in\IP_n$.
\end{enumerate}
\end{proposition}

\begin{proof}
Let us prove first that~\ref{Ideals-1}) holds. Take $b\in\IP_n$.
Suppose that $I_k=\IP_nb\IP_n$. Then due to~\eqref{eq:ranks}, we
obtain that $\rank(b)\geq k$. From the other hand, $b=1\cdot b\cdot
1\in I_k$ and so $\rank(b)\leq k$. Thus, $\rank(b)=k$. Conversely,
suppose that $\rank(b)=k$. Then $b=\bigl(A_i\cup B_i'\bigr)_{1\leq
i\leq k}$ for some partitions $\{1,\ldots,n\}=\bigcup\limits_{1\leq
i\leq k}^{\cdot}A_i$ and $\{1,\ldots,n\}=\bigcup\limits_{1\leq i\leq
k}^{\cdot}B_i$. Take $c\in I_k$ and let $\rank(c)=m\leq k$. Since
\begin{multline}
d=b\tau_{B_1\cup\ldots\cup B_{k+1-m}}=\bigl\{\bigl(A_1\cup\ldots\cup
A_{k+1-m}\bigr)\cup\bigl(B_1\cup\ldots\cup
B_{k+1-m}\bigr)',\\A_{k+2-m}\cup B_{k+2-m}',\ldots,A_k\cup
B_{k}'\bigr\}
\end{multline}
is an element of the rank $m$, then due to~\ref{J=D}) of Theorem
\ref{th:Green}, we obtain that there are $u,v\in\IP_n$ such that
$c=udv=ub\tau_{B_1\cup\ldots\cup B_{k+1-m}}v\in\IP_nb\IP_n$. Thus,
$I_k=\IP_nb\IP_n$ and the proof of~\ref{Ideals-1}) is complete.

Let now $a$ be an arbitrary element of $I$ such that $\rank(a)$ has
the maximum value among the numbers $\rank(x)$, $x\in I$. Then due
to the statement~\ref{Ideals-1}), condition~\eqref{eq:Ideal-Chain}
and the fact that $I=\bigcup\limits_{x\in I}\IP_nx\IP_n$, we have
that $I=I_{\rank(a)}=\IP_na\IP_n$. Thus, statements~\ref{Ideals-2})
and~\ref{Ideals-3}) hold.
\end{proof}

As a corollary we obtain now the following proposition.

\begin{proposition}\label{cor:Ideals}
All the ideals of $\IP_n$ are principal and form the chain
\eqref{eq:Ideal-Chain}.
\end{proposition}

\begin{proof}
Follows from Proposition~\ref{pr:Ideals}.
\end{proof}

Set $\mathcal{D}_k=\bigl\{a\in\IP_n:~\rank(a)=k\bigr\}$ for all
$k\in\mathbb{N}$, $1\leq k\leq n$. Due to~\ref{J=D}) of Theorem
\ref{th:Green}, we have that all these sets exhaust all the
$\GD$-classes of $\IP_n$. Now we are able to formulate a result on
the structure of maximal subsemigroups of $\IP_n$.

\begin{theorem}\label{th:MaximalSusemigroups}
Let $n\geq 3$ and $S$ be a subset of $\IP_n$. Then the following
statements are equivalent:
\begin{enumerate}
\item\label{MaxSubsemi-1}
$S$ is a maximal subsemigroup of $\IP_n$;
\item\label{MaxSubsemi-2}
either $S=\IT_n\cup I_{n-2}$ or $S=G\cup I_{n-1}$ for some maximal
subgroup $G$ of $\Sym_n$.
\end{enumerate}
In addition, every maximal subsemigroup of $\IP_n$ is an inverse
subsemigroup of $\IP_n$.
\end{theorem}

\begin{proof}
Let us prove first that~\ref{MaxSubsemi-2}) implies
\ref{MaxSubsemi-1}). If $S$ coincides with the subsemigroup $G\cup
I_{n-1}$ of $\IP_n$ for some maximal subgroup $G$ of $\Sym_n$ then
since the condition~\eqref{eq:Ideal-Chain}, we have that $S$ is a
maximal subsemigroup of $\IP_n$. Note that $\IT_n\cup I_{n-2}$ is a
subsemigroup of $\IP_n$, as $\IT_n$ is a subsemigroup of $\IP_n$ and
$I_{n-2}$ is an ideal of $\IP_n$. If now $\IT_n\cup I_{n-2}$ is a
proper subsemigroup of $T$, where $T$ is a subsemigroup of $\IP_n$,
then, due to Lemma~\ref{lm:BeforePrGeneratingSystems}, $T$ contains
an element of $\Sym_n\xi_{1,2,3}\Sym_n$ and so, taking to account
Proposition~\ref{pr:Generating-SystemsIP} and the fact that
$\Sym_n\subseteq\IT_n$, we obtain that $T=\IP_n$. Thus,
\ref{MaxSubsemi-2}) implies~\ref{MaxSubsemi-1}).

Let now $S$ be a maximal subsemigroup in $\IP_n$. Note that $S\cup
I_{n-2}$ is a subsemigroup of $\IP_n$. Besides, $S\cup I_{n-2}$ is a
proper subset of $\IP_n$. Indeed, otherwise we would have $S\cup
I_{n-2}=\IP_n$, whence $\Sym_n\cup \mathcal{D}_{n-1}\subseteq S$ and
so due to Proposition~\ref{pr:Generating-SystemsIP}, we would obtain
that $S=\IP_n$. Thus, $S\cup I_{n-2}=S$ and so $I_{n-2}\subseteq S$.
Since $S\cup \{1\}$ is a proper subsemigroup of $\IP_n$, we have
that $S=S\cup\{1\}$ and $G=S\cap \Sym_n\ne\varnothing$. Obviously,
$G$ is a subgroup of $\Sym_n$. Now we have two possibilities.

\emph{Case 1}. $G$ is a proper subgroup of $\Sym_n$. Then
$S\subseteq G\cup I_{n-1}$ and due to the fact that $G\cup I_{n-1}$
is a proper subsemigroup of $\IP_n$, we obtain that $S=G\cup
I_{n-1}$. It remains to note that the latter implies that $G$ is a
maximal subgroup of $\Sym_n$.

\emph{Case 2}. $G=\Sym_n$. Then $\Sym_n\cup I_{n-2}\subseteq S$.
Since $\Sym_n\cup I_{n-2}$ is a proper subsemigroup of $\IT_n\cup
I_{n-2}$, we have that $S$ contains an element $a$ of
$\mathcal{D}_{n-1}$. Then due to
Lemma~\ref{lm:BeforePrGeneratingSystems} and Proposition
\ref{pr:Generating-SystemsIP}, we obtain that $S\subseteq\IT_n\cup
I_{n-2}$. But $\IT_n\cup I_{n-2}$ is a maximal subsemigroup of
$\IP_n$ and so $S=\IT_n\cup I_{n-2}$. This completes the proof of
that~\ref{MaxSubsemi-1}) implies~\ref{MaxSubsemi-2}).

That every maximal subsemigroup of $\IP_n$ is an inverse
subsemigroup of $\IP_n$, follows from what we already have done and
the fact that $\IT_n\cup I_{n-2}$ and $G\cup I_{n-1}$ are inverse
subsemigroups of $\IP_n$ for all subgroups $G$ of $\Sym_n$.
\end{proof}

\section{Automorphism group $\Aut(\IP_X)$}
\label{sec:AutomorphismsIP_X}

Let $g\in\Sym_X$. Denote by $\varphi_g$ the map from $\IP_X$ to
$\IP_X$, given by
\begin{equation}
\varphi_g(a)=g^{-1}ag~\mbox{for every}~a\in \IP_X.
\end{equation}
Clearly, $\varphi_g$ belongs to $\Aut(\IP_X)$, automorphism group of
$\IP_X$. Throughout this section, denote by $\id$ the identity map
of the set $X$ to itself.

The main result of this section is the following theorem.

\begin{theorem}\label{th:Automorphisms}
Let $\varphi\in \Aut(\IP_X)$. Then $\varphi=\varphi_g$ for some
$g\in\Sym_X$. In particular, $\Aut(\IP_X)\cong \Sym_X$ when
$\mid\!X\!\mid\ne 2$ and $\Aut(\IP_2)=\{\id\}$.
\end{theorem}

We will divide the proof of this theorem into few lemmas.

Naturally, $\varphi$ induces an automorphism
$\chi=\varphi\mid_{E(\IP_X)}$ of the semilattice $E(\IP_X)$. Set
$\zeta_{x}=\tau_{X\setminus\{x\}}$ for all $x\in X$. Set also
$\Phi=\bigl\{\zeta_x\in\IP_X:~x\in X\bigr\}$. Recall that if
$(E,\leq)$ is a semilattice with the zero element $0$, then an
element $f$ of $E$ is said to be \emph{primitive} if $g\leq f$
implies either $g=f$ or $g=0$, for all $g\in E$. For all $n\geq 2$
set
\begin{multline}
\Theta_{\max}^n=\bigl\{\tau_{i,j}\in\IP_n:~i,j\in\{1,\ldots,n\},~i\ne
j\bigr\}~\mbox{and}\\\Theta_{\mathrm{pr}}^n=\bigl\{\tau_{F}\tau_{\{1,\ldots,n\}\setminus
F}\in\IP_n:F~\mbox{is a proper subset
of}~\{1,\ldots,n\}\bigr\}=\\\mathcal{D}_2\cap E(\IP_n).
\end{multline}
Notice that $\Phi\subseteq\Theta_{\mathrm{pr}}^n$.

\begin{lemma}\label{lm:Aut-Max-and-Pr}
Let $n\geq 2$. Then the set of all primitive elements of the
semilattice $E(\IP_n)$ coincides with $\Theta_{\mathrm{pr}}^n$. Also
then the set of all maximal elements of the semilattice
$E(\IP_n)\setminus\{1\}$ coincides with $\Theta_{\max}^n$.
\end{lemma}

\begin{proof}
Follows from Proposition~\ref{pr:Omega}.
\end{proof}

\begin{lemma}\label{lm:Aut-Semilattice}
Take $\theta\in\Aut\bigl(E(\IP_n)\bigr)$. Then there is $g\in\Sym_n$
such that $\theta(e)=g^{-1}eg$ for all $e\in E(\IP_n)$.
\end{lemma}

\begin{proof}
Clearly, the statement holds when $n=1$. Thus, let us assume that
$n\geq 2$.

Obviously, $\theta(1)=1$. Then
$\theta\bigl(E(\IP_n)\setminus\{1\}\bigr)=E(\IP_n)\setminus\{1\}$.
Hence, due to Lemma~\ref{lm:Aut-Max-and-Pr}, we obtain that
$\theta\bigl(\Theta_{\max}^n\bigr)=\Theta_{\max}^n$ and
$\theta\bigl(\Theta_{\mathrm{pr}}^n\bigr)=\Theta_{\mathrm{pr}}^n$.
Take $f=\tau_{F}\tau_{\{1,\ldots,n\}\setminus
F}\in\Theta_{\mathrm{pr}}^n$. Set
$\Lambda_f=\bigl\{a\in\Theta_{\max}^n:~fa=f\bigr\}$. Then
$\theta(\Lambda_f)=\Lambda_{\theta(f)}$. If $f\notin\Phi$ then
$2\leq\mid\!F\!\mid\leq n-2$. Thus,
\begin{equation}\label{eq:Lambda_f}
\mid\!\Lambda_f\!\mid=\binom{\mid\!F\!\mid}{2}+\binom{n-\mid\!F\!\mid}{2},
~\mbox{if}~f\notin\Phi.
\end{equation}
Otherwise, we have the following:
\begin{equation}\label{eq:Zeta}
\mid\!\Lambda_f\!\mid=\binom{n-1}{2},~\mbox{if}~f\in\Phi.
\end{equation}
Let us prove now that for all $n\geq 4$ and for all $k$, $2\leq
k\leq n-2$, the following holds:
\begin{equation}\label{eq:Inequality}
\binom{k}{2}+\binom{n-k}{2}<\binom{n-1}{2}.
\end{equation}
Indeed, the inequality
\begin{multline}
k(k-n)=(k^2-1)+1-kn=(k-1)(k+1)+1-kn<\\(k-1)n+1-kn=1-n~\mbox{implies}
\end{multline}
\begin{multline}
\binom{k}{2}+\binom{n-k}{2}=\frac{1}{2}\bigl(k(k-1)+(n-k)(n-k-1)\bigr)=\\\frac{1}{2}
\bigl(2k^2-2kn+n^2-n\bigr)=k(k-n)+\frac{1}{2}
\bigl(n^2-n\bigr)<\\\frac{1}{2}
\bigl(n^2-n\bigr)+1-n=\frac{1}{2}(n-1)(n-2)=\binom{n-1}{2}.
\end{multline}
Now due to~\eqref{eq:Lambda_f},~\eqref{eq:Zeta},
\eqref{eq:Inequality} and the equality
$\theta(\Lambda_f)=\Lambda_{\theta(f)}$, we obtain that
$\theta(\Phi)=\Phi$. Then there is an element $g$ of $\Sym_n$ such
that $\theta(\zeta_{x})=\zeta_{g(x)}$ for all $x\in\{1,\ldots,n\}$.

Take now distinct $x$ and $y$ of $\{1,\ldots,n\}$. Since
$\zeta_{x}\tau_{x,y}=0$ and $\zeta_{y}\tau_{x,y}=0$, we have that
$\zeta_{g(x)}\theta(\tau_{x,y})=0$ and
$\zeta_{g(y)}\theta(\tau_{x,y})=0$. The latter, taking to account
$\theta\bigl(\Theta_{\max}^n\bigr)=\Theta_{\max}^n$, implies that
$\theta(\tau_{x,y})=\tau_{g(x),g(y)}=g^{-1}\tau_{x,y}g$.

Let now $e=(E_i\cup E_i')_{i\in I}$ be a nonidentity idempotent
element of $E(\IP_n)$. Then
\begin{equation}
e=\prod\bigl\{\tau_{x,y}:~x\ne y,~\{x,y\}\subseteq E_i~\mbox{for
some}~i\in I\bigr\}
\end{equation}
implies
\begin{multline}
\theta(e)=\prod\bigl\{\tau_{g(x),g(y)}:~x\ne y,~\{x,y\}\subseteq
E_i~\mbox{for some}~i\in
I\bigr\}=\\\prod\bigl\{g^{-1}\tau_{x,y}g:~x\ne y,~\{x,y\}\subseteq
E_i~\mbox{for some}~i\in I\bigr\}=g^{-1}eg.
\end{multline}
This completes the proof.
\end{proof}

Take distinct $x$ and $y$ of $X$. Define an element
$\varepsilon_{x,y}$ of $\Sym_X$ as follows:
\begin{equation}
\varepsilon_{x,y}(x)=y,~\varepsilon_{x,y}(y)=x~\mbox{and}~\varepsilon_{x,y}(t)=t~\mbox{for
all}~t\in X\setminus\{x,y\}.
\end{equation}

\begin{corollary}\label{cor:Aut-IP_6}
Let $\mid\!X\!\mid=6$. Then there is $g\in\Sym_6$ such that
$\varphi(h)=\varphi_g(h)$ for all $h\in\Sym_6$.
\end{corollary}

\begin{proof}
If we put $\chi=\theta$ and $n=6$ in the statement of Lemma
\ref{lm:Aut-Semilattice}, we will obtain that there is $g\in\Sym_6$
such that $\chi(e)=g^{-1}eg$ for all $e\in E(\IP_6)$. Take distinct
$x$ and $y$ of $\{1,\ldots,6\}$. Then
\begin{equation}
g^{-1}\tau_{x,y}g=
\varphi(\tau_{x,y})=\varphi(\tau_{x,y}\varepsilon_{x,y})=g^{-1}\tau_{x,y}g\varphi(\varepsilon_{x,y}),
\end{equation}
whence
\begin{equation}
\tau_{x,y}=\tau_{x,y}g\varphi(\varepsilon_{x,y})g^{-1}.
\end{equation}
The latter implies that either $g\varphi(\varepsilon_{x,y})g^{-1}=1$
or $g\varphi(\varepsilon_{x,y})g^{-1}=\varepsilon_{x,y}$. But since
the order of $g\varphi(\varepsilon_{x,y})g^{-1}$ equals $2$, we have
that $g\varphi(\varepsilon_{x,y})g^{-1}=\varepsilon_{x,y}$, which is
equivalent to $\varphi(\varepsilon_{x,y})=g^{-1}\varepsilon_{x,y}g$.
Now, taking to account the known fact that
$\langle\varepsilon_{x,y}:~x\ne y\rangle=\Sym_n$ (see
\cite{Kurosh}), we obtain that $\varphi(h)=\varphi_g(h)$ for all
$h\in\Sym_6$.
\end{proof}

\begin{lemma}\label{lm:Aut-S_X}
There is $g\in\Sym_X$ such that $\varphi(h)=\varphi_g(h)$ for all
$h\in\Sym_X$.
\end{lemma}

\begin{proof}
Due to Corollary~\ref{cor:Aut-IP_6}, we have that the statement
holds when $\mid\!X\!\mid=6$. Assume now that $\mid\!X\!\mid\ne 6$.

Since $\varphi$ preserves the set of all invertible elements of
$\IP_X$, we have, due to Proposition~\ref{cor:Invertible-Elements},
that $\varphi(\Sym_X)=\Sym_X$. Hence, $\varphi$ induces an
automorphism of $\Sym_X$. Then due to known fact, which claims that
if $\mid\!X\!\mid\ne 6$ then every automorphism of $\Sym_X$ is inner
(see \cite{Kurosh}), we have that there is $g\in\Sym_X$ such that
$\varphi(h)=g^{-1}hg=\varphi_g(h)$ for all $h\in\Sym_X$. This
completes the proof.
\end{proof}

Set now $\psi =\varphi\varphi_g$. Then $\psi$ is, obviously, an
automorphism of $\IP_X$ and, due to Lemma~\ref{lm:Aut-S_X},
$\psi\mid_{\Sym_X}$ is the identity map of $\Sym_X$ to itself. For
all $M\subseteq X$ set
\begin{equation}
\widetilde{\Sym}_M=\bigl\{h\in\Sym_X:~h(x)=x~\mbox{for all}~x\in
X\setminus M\bigr\}.
\end{equation}
For all $a\in\IP_X$ set
\begin{equation}
\Fix_l(a)=\bigl\{h\in\Sym_X:~ha=a\bigr\}~\mbox{and}~\Fix_r(a)=\bigl\{h\in\Sym_X:~ah=a\bigr\}.
\end{equation}

\begin{lemma}\label{lm:Aut-C(e)}
Let $a\in\IP_X$. Let also $X=\bigcup\limits_{i\in
I}^{\cdot}A_i=\bigcup\limits_{i\in I}^{\cdot}B_i$. Then
\begin{enumerate}
\item
$\Fix_l(a)=\bigoplus\limits_{i\in I}\widetilde{\Sym}_{A_i}$ if and
only if $a=\bigl(A_i\cup U_i'\bigr)_{i\in I}$ for some partition
$X=\bigcup\limits_{i\in I}^{\cdot}U_i$;
\item
$\Fix_r(a)=\bigoplus\limits_{i\in I}\widetilde{\Sym}_{B_i}$ if and
only if $a=\bigl(V_i\cup B_i'\bigr)_{i\in I}$ for some partition
$X=\bigcup\limits_{i\in I}^{\cdot}V_i$.
\end{enumerate}
\end{lemma}

\begin{proof}
Straightforward.
\end{proof}

\begin{corollary}\label{cor:Aut-H-classes}
$a\GH\psi(a)$ for all $a\in\IP_X$. In particular, $\psi(e)=e$ for
all $e\in E(\IP_X)$.
\end{corollary}

\begin{proof}
That $a\GH\psi(a)$ for all $a\in\IP_X$ follows from Lemma
\ref{lm:Aut-C(e)} and Theorem~\ref{th:Green}. Then $\psi(e)=e$ for
all $e\in E(\IP_X)$, due to the fact that every $\GH$-class of an
arbitrary semigroup contains at most one idempotent (see Corollary
2.2.6 from \cite{Howie}).
\end{proof}

\begin{lemma}\label{lm:Aut-n-geq-3}
Let $a\in\IP_X$ and $\rank(a)\geq 3$. Then $\psi(a)=a$.
\end{lemma}

\begin{proof}
Let $a=\bigl(A_i\cup B_i'\bigr)_{i\in I}$, $\mid\!I\!\mid\geq 3$.
Due to Corollary~\ref{cor:Aut-H-classes}, we have that $a\GH\psi(a)$
and so $\psi(a)=\bigl(A_i\cup B_{\alpha(i)}'\bigr)_{i\in I}$ for
some bijective map $\alpha:I\to I$. Due to Corollary
\ref{cor:Aut-H-classes}, we also have that $ea\GH e\psi(a)$ for all
$e\in E(\IP_X)$.

Take arbitrary distinct $i$ and $j$ of $I$. Since $\tau_{A_i\cup
A_j}a\GH\tau_{A_i\cup A_j}\psi(a)$, we have that
\begin{multline}
\bigl\{\bigl(A_i\cup A_j\bigr)\cup\bigl(B_i\cup
B_j\bigr)',\bigl(A_l\cup B_l'\bigr)_{l\in
I\setminus\{i,j\}}\bigr\}~\mbox{and}\\\bigl\{\bigl(A_i\cup
A_j\bigr)\cup\bigl(B_{\alpha(i)}\cup
B_{\alpha(j)}\bigr)',\bigl(A_l\cup B_l'\bigr)_{l\in
I\setminus\{i,j\}}\bigr\}
\end{multline}
are $\GH$-equivalent, whence $\{i,j\}=\{\alpha(i),\alpha(j)\}$. Let
now $k\in I$. Then $\alpha(k)=k$. Suppose the contrary. Then
$\{k,m\}=\{\alpha(k),\alpha(m)\}$ for all $m\in I\setminus\{k\}$
implies that $\alpha(k)=m$ for all $m\in I\setminus\{k\}$. But
$\mid\!I\!\mid\geq 3$ and we get a contradiction. Thus, $\alpha$ is
an identity map of $I$ to itself, which is equivalent to
$\psi(a)=a$. This completes the proof.
\end{proof}

Note that since $\IP_1$ is isomorphic to the unit group and since
$\IP_2\cong\mathbb{Z}_2^0$, where $\mathbb{Z}_2^0$ denotes the group
$\mathbb{Z}_2$ with adjoint zero, we have that $\Aut(\IP_X)=\{\id\}$
when $\mid\!X\!\mid\leq 2$.

\begin{lemma}\label{lm:Aut-leq-2}
Let $a\in\IP_X$ and $\rank(a)\leq 2$. Then $\psi(a)=a$.
\end{lemma}

\begin{proof}
If $\rank(a)=1$ then $a=0$ and, obviously, $\psi(a)=a$. So let us
suppose that $\rank(a)=2$. Assume that $a=\bigl\{A\cup B',C\cup
D'\bigr\}$. Fix $x\in A$ and $y\in B$.

Suppose that $\lmod A\rmod\geq 2$ and $\lmod B\rmod\geq 2$. Then
$\psi(a)=a$. Indeed, we have that $A\setminus\{x\}\ne\varnothing$
and $B\setminus\{y\}\ne\varnothing$, so if $y_1\in B\setminus\{y\}$
then we can consider the equality
$a=\bigl\{\{x,y'\},\bigl(A\setminus\{x\}\bigr)\bigcup\bigl(B\setminus\{y\}\bigr)',C\cup
D'\bigr\}\cdot\tau_{y,y_1}$, whence, due to Corollary
\ref{cor:Aut-H-classes} and Lemma~\ref{lm:Aut-n-geq-3}, we will have
that $\psi(a)=a$.

Analogously, if $\lmod C\rmod\geq 2$ and $\lmod D\rmod\geq 2$ then
$\psi(a)=a$.

Thus, we may assume that either $\lmod A\rmod=1$ or $\lmod
B\rmod=1$, and that either $\lmod C\rmod=1$ or $\lmod D\rmod=1$.
Without loss of generality we may suppose that $\lmod A\rmod=1$.
Then we will have two possibilities.

\emph{Case 1}. $\lmod C\rmod=1$. Then $\lmod X\rmod=2$ and we obtain
$\psi=\id$.

\emph{Case 2}. $\lmod D\rmod=1$. Then $\psi(a)=a$. Suppose the
contrary. Then we would obtain that $\psi(a)=\bigl\{A\cup D',C\cup
B'\bigr\}=\zeta_{x}h$ for some $h\in\Sym_X$. But
$\zeta_{x}h=\psi(\zeta_{x}h)$ and so $a=\zeta_{x}h$, whence $B=D$,
which leads to a contradiction.

Thus, we proved that $\psi(a)=a$, which was required.
\end{proof}

As a consequence of that we have from Lemmas~\ref{lm:Aut-n-geq-3}
and~\ref{lm:Aut-leq-2}, we have that $\psi=\id$, whence
$\varphi=\varphi_g$. It remains to prove that
$\Aut(\IP_X)\cong\Sym_X$ when $\lmod X\rmod\ne 2$. This follows from
the following lemma.

\begin{lemma}\label{lm:Aut-cong-Sym_X}
Suppose that $\lmod X\rmod\geq 3$. Then a map
$\vartheta:~\Sym_X\to\Aut(\IP_X)$, given by
\begin{equation}
\vartheta(h)=\varphi_h~\mbox{for all}~h\in\Sym_X,
\end{equation}
is an isomorphism from $\Sym_X$ onto $\Aut(\IP_X)$.
\end{lemma}

\begin{proof}
We have already proved that $\vartheta$ is an onto homomorphism from
$\Sym_X$ to $\Aut(\IP_X)$. But, besides, $\vartheta$ is an injective
map. Indeed, $\vartheta(h_1)=\vartheta(h_2)$ implies that
$h_1^{-1}hh_1=h_2^{-1}hh_2$ or just that
$\bigl(h_1h_2^{-1}\bigr)^{-1}h\bigl(h_1h_2^{-1}\bigr)=h$ for all
$h\in\Sym_X$ and it remains to note that $\Sym_X$ is a center-free
group when $\lmod X\rmod\geq 3$ (see \cite{Kurosh}). Thus,
$\vartheta$ is an isomorphism.
\end{proof}

The proof of theorem is complete.

\section{Connections between $\IP_X$ and other semigroups}\label{sec:Connections}

Set $\Upsilon=\bigl\{X,X'\bigr\}$. Then $\Upsilon\in\Ch_X$. The
following proposition shows that $\IP_X\cup\{\Upsilon\}$ is a
maximal inverse subsemigroup of $\Ch_X$ when $\lmod X\rmod\geq 2$.

\begin{proposition}\label{pr:IP_X-inverse-subsemigroup}
Let $\lmod X\rmod\geq 2$. Then $\IP_X\cup \{\Upsilon\}$ is a maximal
inverse subsemigroup of $\Ch_X$.
\end{proposition}

\begin{proof}
Since $\IP_X$ is an inverse subsemigroup of $\Ch_X$ and
$a\Upsilon=\Upsilon a=\Upsilon$ for all $a\in\IP_X\cup
\{\Upsilon\}$, we obtain that $\IP_X\cup \{\Upsilon\}$ is a proper
inverse subsemigroup of $\Ch_X$.

Suppose now that $S$ is an inverse subsemigroup of $\Ch_X$ such that
$\IP_X\cup\{\Upsilon\}$ is a subsemigroup of $S$. Take $s\in
S\setminus\IP_X$. Then there is a nonempty subset $A$ of $X$ such
that either $s$ contains a block $A$ or $s$ contains a block $A'$.
Without loss of generality we may assume that $s$ contains the block
$A$. Let $t$ be the inverse of $s$ in $S$. Then $st$ is an
idempotent in $S$ and so, due to the fact that idempotents of
inverse semigroup commute, we obtain that
$u=st\cdot\Upsilon=\Upsilon\cdot st$. The latter implies that $u$
contains both blocks $A$ and $X$, whence $A=X$. Then $s$ is an
idempotent and due to equalities $s=\Upsilon s$ and $\Upsilon
s=s\Upsilon$, we have that $s$ contains the block $X'$ and so
$s=\Upsilon$. That is, $S=\IP_X\cup\{\Upsilon\}$. This implies that
$\IP_X\cup\{\Upsilon\}$ is a maximal inverse subsemigroup of $\Ch_X$
which was required.
\end{proof}

Denote by $\IS_X$ the \emph{symmetric inverse} semigroup on the set
$X$. Let $s\in\IS_X$. Denote by $\dom(s)$ and $\ran(s)$ the
\emph{domain} and the \emph{range} of $s$ respectively. The
following theorem shows how one can embed the symmetric inverse
semigroup into the inverse partition one.

\begin{theorem}\label{th:IS-embeds-into-IP}
Let $\overline{x}\notin X$. Then $\IS_X$ isomorphically embeds into
$\IP_{X\cup\{\overline{x}\}}$.
\end{theorem}

\begin{proof}
For all $s\in\IS_X$, set
\begin{equation}\label{eq:Omega}
\Omega_s=\bigl(X\cup\{\overline{x}\}\setminus\dom(s)\bigr)\bigcup
\bigl(X\cup\{\overline{x}\}\setminus\ran(s)\bigr)'.
\end{equation}
Set a map $\kappa:~\IS_X\to\IP_{X\cup\{\overline{x}\}}$ as follows:
\begin{equation}\label{eq:kappa-definition}
\kappa(s)=\bigl\{\Omega_s,\bigl(\{x,s(x)'\}\bigr)_{x\in\dom(s)}\bigr\}~\mbox{for
all}~s\in\IS_X.
\end{equation}
Take an arbitrary $s$ of $\IS_X$. Then we have the following
condition:
\begin{equation}\label{eq:kappa}
x\equiv_{\kappa(s)}\overline{x}\equiv_{\kappa(s)}{\overline{x}}'\equiv_{\kappa(s)}y'~\mbox{for
all}~x\in X\setminus\dom(s)~\mbox{and}~y\in X\setminus\ran(s).
\end{equation}
Take $s,t\in\IS_X$. Then due to~\eqref{eq:Omega} and
\eqref{eq:kappa}, we obtain that
\begin{multline}\label{eq:dreadful}
x\equiv_{\kappa(s)\kappa(t)}\overline{x}\equiv_{\kappa(s)\kappa(t)}\overline{x}'
\equiv_{\kappa(s)\kappa(t)}y'~\mbox{for all}~x,y\in X~\mbox{such
that}\\x\notin
s^{-1}\bigl(\ran(s)\cap\dom(t)\bigr)~\mbox{and}~y\notin
t\bigl(\dom(t)\cap\ran(s)\bigr).
\end{multline}
Notice that
\begin{equation}\label{eq:S-T}
s^{-1}\bigl(\ran(s)\cap\dom(t)\bigr)=\dom(st)~\mbox{and}~
t\bigl(\dom(t)\cap\ran(s)\bigr)=\ran(st).
\end{equation}
If now $x\in\dom(st)$ then
$x\equiv_{\kappa(s)}s(x)'~\mbox{and}~s(x)\equiv_{\kappa(t)}st(x)'$,
whence
\begin{equation}\label{eq:last}
x\equiv_{\kappa(s)\kappa(t)}st(x)'~\mbox{for all}~x\in\dom(st).
\end{equation}
The conditions~\eqref{eq:dreadful},~\eqref{eq:S-T} and
\eqref{eq:last} imply that
\begin{equation}
\kappa(s)\kappa(t)=\bigl\{\Omega_{st},\bigl(\{x,
st(x)'\}\bigr)_{x\in\dom(st)}\bigr\}=\kappa(st).
\end{equation}
Thus, $\kappa$ is a homomorphism from $\IS_X$ to
$\IP_{X\cup\{\overline{x}\}}$. It remains to prove that $\kappa$ is
an injective map.

Suppose that $\kappa(s)=\kappa(t)$ for some $s,t\in\IS_X$. Then it
follows from~\eqref{eq:kappa-definition} that
$\dom(s)\subseteq\dom(t)$ and $\dom(t)\subseteq\dom(s)$, whence
$\dom(s)=\dom(t)$. Then~\eqref{eq:kappa-definition} implies that
$s(x)=t(x)$ for all $x\in\dom(s)=\dom(t)$. Hence, $s=t$ and so
$\kappa$ is injective. The proof is complete.
\end{proof}

It follows immediately from Theorem~\ref{th:IS-embeds-into-IP} that
$\IS_n$ embeds into $\IP_{n+1}$ for all $n\in\mathbb{N}$.
Surprisingly, the following theorem shows that one can not construct
an embedding map from $\IS_n$ to $\IP_n$.

\begin{theorem}\label{th:IS-doesn't-embeds-into-IP}
Let $n\in\mathbb{N}$. There is no an injective homomorphism from
$\IS_n$ to $\IP_n$.
\end{theorem}

\begin{proof}
Suppose the contrary. Then there is a subsemigroup $U$ of $\IP_n$
such that $U\cong\IS_n$. Then we have that $U$ is a regular
subsemigroup of $\IP_n$, whence, due to Proposition 2.4.2 from
\cite{Howie}, we obtain that $\GD^{U}=\GD\cap(U\times U)$, where
$\GD^{U}$ denotes the Green's $\GD$-relation on $U$. Note that
$\IP_n$ contains exactly $n$ different $\GD$-classes. This implies
that $U$ contains at most $n$ different $\GD^{U}$-classes. But since
$U\cong\IS_n$, we have that $U$ contains exactly $n+1$ different
$\GD^{U}$-classes. We get a contradiction. This completes the proof.
\end{proof}

\section{$\IP_n$ embeds into $\IS_{2^n-2}$}\label{ef}

Let $S$ be an inverse semigroup with the natural partial order
$\leq$ on it. For $A\subseteq S$ denote by $[A]$ the order ideal of
$S$ with respect to $\leq$, i.e., $[A]=\bigl\{b:~a\leq b~\mbox{for
some}~a\in A\bigr\}$. Let also $H$ be a \emph{closed inverse
subsemigroup} of $S$, i.e., $H$ is an inverse subsemigroup of $S$
and $[H]=H$ (see \cite{Howie}). Recall (see \cite{Howie}) that one
can define the set of all \emph{right $\leq$-cosets} of $H$ as
follows:
\begin{equation}
\mathcal{C}=\mathcal{C}_H=\bigl\{[Hs]:~ss^{-1}\in H\bigr\}.
\end{equation}
Further, one can define the \emph{effective transitive
representation} $\phi_H:S\to\IS_{\mathcal{C}}$, given by
\begin{equation}
\phi_H(s)=\bigr\{\bigl([Hx],[Hxs]\bigr):~
[Hx],[Hxs]\in\mathcal{C}\bigl\}.
\end{equation}
Let now $K$ and $H$ be arbitrary closed inverse subsemigroups of
$S$. For a definition of the \emph{equivalence} of representations
$\phi_K$ and $\phi_H$, we refer reader to \cite{Howie}. But we note
that due to Proposition IV.4.13 from \cite{Petrich}, one has that
$\phi_K$ and $\phi_H$ are equivalent if and only if there exists
$a\in S$ such that $a^{-1}Ha\subseteq K$ and $aKa^{-1}\subseteq H$.
We will need the following well-known fact.

\begin{theorem}[Proposition 5.8.3 from \cite{Howie}]\label{th:Howie}
Let $H$ be a closed inverse subsemigroup of an inverse semigroup $S$
and let $a,b\in S$. Then $[Ha]=[Hb]$ if and only if $ab^{-1}\in H$.
\end{theorem}

The main result of this section is the following theorem.

\begin{theorem}
Let $n\geq 2$. Up to equivalence, there is only one faithful
effective transitive representation of $\IP_n$, namely to
$\IS_{2^n-2}$. In particular, $\IP_n$ isomorphically embeds into
$\IS_{2^n-2}$.
\end{theorem}

We divide the proof of this theorem into lemmas. Throughout all
further text of this section we suppose that $H$ is a closed inverse
subsemigroup of $\IP_n$.

\begin{lemma}\label{lm:Ef-Tr-1}
$H=[G]$ for some subgroup $G$ of $\IP_n$.
\end{lemma}

\begin{proof}
Since $\IP_n$ is finite, we have that $E(H)$ contains a zero
element. It remains to use Proposition IV.5.5 from \cite{Petrich},
which claims that if the set of idempotents of a closed inverse
subsemigroup contains a zero element, then this subsemigroup is a
closure of some subgroup of the original semigroup.
\end{proof}

Denote by $e$ the identity element of $G$.

\begin{lemma}\label{lm:Ef-Tr-2}
If $e=0$ then $\phi_H$ is not faithful.
\end{lemma}

\begin{proof}
We have $G=\{0\}$, whence $H=[0]=\IP_n$ and so $[Hx]\supseteq
[0]=\IP_n$ for all $x\in\IP_n$. Thus, $[Hx]=\IP_n$ for all
$x\in\IP_n$. Then $\lmod\phi_H(\IP_n)\rmod=1$, whence we obtain that
$\phi_H$ is not faithful.
\end{proof}

\begin{lemma}\label{lm:Ef-Tr-3}
Let $\rank(e)\geq 3$. Then $\phi_H$ is not faithful.
\end{lemma}

\begin{proof}
Take $b\in\mathcal{D}_2$. Since $bb^{-1}\in\mathcal{D}_2$, we have
that $bb^{-1}\notin H$ and so $[Hb]\notin\mathcal{C}$. The latter
gives us that $\phi_H(b)$ equals the zero element of
$\IS_{\mathcal{C}}$. Then, due to $\lmod\mathcal{D}_2\rmod\geq 2$,
we obtain that $\phi_H$ is not faithful.
\end{proof}

\begin{lemma}\label{lm:Ef-Tr-4}
Let $\rank(e)=2$ and $G\cong\mathbb{Z}_2$. Then $\phi_H$ is not
faithful.
\end{lemma}

\begin{proof}
Let $G=\{e,q\}$. We are going to prove that $\phi_H(e)=\phi_H(q)$.

Let us prove first that
$\dom\bigl(\phi_H(e)\bigr)=\dom\bigl(\phi_H(q)\bigr)$. Indeed, take
$[Hx]\in\mathcal{C}$. Then, due to the equality
$(xe)(xe)^{-1}=xex^{-1}=xqq^{-1}x^{-1}=(xq)(xq)^{-1}$, we obtain
that $[Hxe]\in\mathcal{C}$ if and only if $(xe)(xe)^{-1}\in H$ if
and only if $(xq)(xq)^{-1}\in H$ if and only if
$[Hxq]\in\mathcal{C}$. Thus,
$\dom\bigl(\phi_H(e)\bigr)=\dom\bigl(\phi_H(q)\bigr)$.

Take now $x\in\dom\bigl(\phi_H(e)\bigr)$. Then $xex^{-1}\in
H=[\{e,q\}]$. But since $xex^{-1}$ is an idempotent and
$\rank(xex^{-1})\leq\rank(e)=2$, we obtain, taking to account
Proposition~\ref{pr:Omega}, that $xex^{-1}=e$. Hence,
$(xe)(xe)^{-1}=ee^{-1}$ and so, due to Proposition~2.4.1 from
\cite{Howie}, we obtain that $xe\GR e$. But then we have that
$\rank(xe)=\rank(e)$ and due to $\lambda_{xe}\supseteq\lambda_{e}$
(which follows, in turn, from \eqref{eq:remark-for-ranks}), we
deduce that $\lambda_{xe}=\lambda_{e}$, whence due to
Theorem~\ref{th:Green}, we have that $xe\GL e$. Thus, $xe\GH e$,
whence $xe\in G$ and so $xq=xe\cdot q\in G$. But then
$(xq)(xe)^{-1}\in G\subseteq H$, whence, due to
Theorem~\ref{th:Howie}, we have that $[Hxe]=[Hxq]$. The latter
implies that $\phi_H(e)(x)=\phi_H(q)(x)$. Thus,
$\phi_H(e)=\phi_H(q)$ and so $\phi$ is not faithful.
\end{proof}

\begin{lemma}\label{ugly-lemma}
Let $f\in\Theta_{\mathrm{pr}}^{n}$ and $T=[f]$. Take
$[Tx]\in\mathcal{C}_T$. Then $\rank(fx)=2$ and $[Tx]=[fx]$.
\end{lemma}

\begin{proof}
Clearly, $[Tx]\in\mathcal{C}_T$ is equivalent to $f\leq xx^{-1}$.

Obviously, $\rank(fx)\leq\rank(f)=2$. But $\rank(fx)=1$ is
impossible. Indeed, otherwise we would have $fx=0$, whence
$0=fxx^{-1}=f$, which does not hold. Thus, $\rank(fx)=2$.

Note that $[fx]\subseteq[Tx]$. It remains to prove that
$[Tx]\subseteq[fx]$. Take $t\in T$. Then $f\leq t$ and due to the
fact that the natural partial order on an arbitrary inverse
semigroup is compatible (see \cite{Howie}), we obtain that $fx\leq
tx$. That is, $tx\in [fx]$. Hence, $Tx\subseteq[fx]$, whence
$[Tx]\subseteq\bigl[[fx]\bigr]=[fx]$.

The proof is complete.
\end{proof}

\begin{lemma}\label{lm:Ef-Tr-5}
Let $\rank(e)=2$ and $G=\{e\}$. Then $\phi_H$ is faithful.
\end{lemma}

\begin{proof}
Note that $H=[e]$. Let $e=\tau_E\tau_{E_1}$, where $E$ and $E_1$ are
nonempty subsets of $\{1,\ldots,n\}$ such that
$\{1,\ldots,n\}=E\bigcup\limits^{\cdot}E_1$. Suppose that
$\phi_H(s)=\phi_H(t)$ for some $s$ and $t$ of $\IP_n$. Let $A$ be an
arbitrary $\rho_{t}$-class. Set
$\overline{A}=\{1,\ldots,n\}\setminus A$.

Suppose first that $s=0$. We are going to prove that $t=0$. Suppose
the contrary. We have that $\rank(e\cdot xs)=1$ for all $x\in\IP_n$
such that $xx^{-1}\in [e]$. So, due to Lemma~\ref{ugly-lemma}, we
obtain that $\dom\bigl(\phi_H(s)\bigr)=\varnothing$. Then, again by
Lemma~\ref{ugly-lemma}, we have that $\rank(e\cdot xt)=1$, or just
that $ext=0$, for all $x\in\IP_n$ such that $xx^{-1}\in [e]$. Put
now $u=\bigl\{E\cup A',E_1\cup\overline{A}'\bigr\}$ (note that, due
to assumption, $\overline{A}\ne\varnothing$). Then $uu^{-1}=e\in
[e]$ and $eut\ne 0$. Thus, we get a contradiction and so $s=0$
implies $t=0$. Analogously, $t=0$ implies $s=0$.

Assume now that $s\ne 0$, then $t\ne 0$ and so
$\overline{A}\ne\varnothing$. Put again $u=\bigl\{E\cup
A',E_1\cup\overline{A}'\bigr\}$. Due to Theorem~\ref{th:Howie} and
the equality $\phi_H(s)=\phi_H(t)$, we have that $(xt)(xs)^{-1}\in
H$ for all $x\in\dom\bigl(\phi_H(t)\bigr)$. Note that
$u\in\dom\bigl(\phi_H(t)\bigr)$. Indeed, we have $uu^{-1}=e\in [e]$
and since $A$ is a $\rho_{tt^{-1}}$-class, we have that
\begin{equation}
(ut)(ut)^{-1}=utt^{-1}u^{-1}=e\in [e].
\end{equation}
This implies that $u\cdot ts^{-1}\cdot u^{-1}\in [e]$. Moreover,
since $\rank(uts^{-1}u^{-1})\leq\rank(u)=2$, we obtain that
$\rank(uts^{-1}u^{-1})=2$, whence $(ut)(us)^{-1}=uts^{-1}u^{-1}=e$.
In particular, we have that $us\ne 0$. But then $A$ is a union of
some $\rho_s$-classes. Since $A$ was an arbitrary chosen
$\rho_t$-class, we obtain that $\rho_s\subseteq\rho_t$. Analogously,
one can prove that $\rho_t\subseteq\rho_s$. Thus, $\rho_s=\rho_t$.
Further, if $s$ contains a block $A\cup B'$ then $us=\bigl\{E\cup
B',E_1\cup \overline{B}'\bigr\}$, where
$\overline{B}=\{1,\ldots,n\}\setminus B$. But $ut=\bigl\{E\cup
A',E_1\cup \overline{A}'\bigr\}$ and so
\begin{multline}
\bigl\{E\cup E',E_1\cup E_1'\bigr\}=e=(ut)(us)^{-1}=\bigl\{E\cup
A',E_1\cup \overline{A}'\bigr\}\cdot\bigl\{E\cup B',E_1\cup
\overline{B}'\bigr\}^{-1}=\\\bigl\{E\cup A',E_1\cup
\overline{A}'\bigr\}\cdot\bigl\{B\cup E',\overline{B}\cup
E_1'\bigr\}.
\end{multline}
This implies $A=B$. Indeed, otherwise we would have
$B\subseteq\overline{A}$ and so $A\subseteq\overline{B}$, whence
$e=\bigl\{E\cup E_1',E_1\cup E'\bigr\}$, which is not true. Again,
since $A$ was an arbitrary chosen $\rho_t$-class, we have that
$\equiv_s=\equiv_t$. Thus, $s=t$. The proof is complete.
\end{proof}

\begin{lemma}\label{lm:Ef-Tr-6}
Let $f\in\Theta_{\mathrm{pr}}^{n}$. Then
$\lmod\mathcal{C}_{[f]}\rmod=2^n-2$.
\end{lemma}

\begin{proof}
Take $[Hx]$ and $[Hy]$ of $\mathcal{C}_{[f]}$. Then due to
Lemma~\ref{ugly-lemma}, we have that $[fx]=[fy]$ and
$\rank(fx)=\rank(fy)$, whence $fx=fy$. Conversely, if $fx=fy$ then
$[Hx]=[fx]=[fy]=[Hy]$. Thus, since $\rank(fx)=\rank(f)$ and $fx=f$
hold simultaneously if and only if $f\GL fx$, we obtain that
$\lmod\mathcal{C}_{[f]}\rmod$ equals the cardinality of $\GL$-class,
which contains $f$, which, in turn, equals the number of all
partitions of $\{1,\ldots,n\}$ into two nonempty blocks. The latter
number is equal to $2^n-2$.
\end{proof}

\begin{lemma}\label{lm:Ef-Tr-7}
Let $f_1,f_2\in\Theta_{\mathrm{pr}}^{n}$. Then $\phi_{[f_1]}$ and
$\phi_{[f_2]}$ are equivalent.
\end{lemma}

\begin{proof}
Let $f_1=\tau_{F_1}\tau_{\{1,\ldots,n\}\setminus F_1}$ and
$f_2=\tau_{F_2}\tau_{\{1,\ldots,n\}\setminus F_2}$ for certain
proper subsets $F_1$ and $F_2$ of $\{1,\ldots,n\}$. Put
$a=\bigl\{F_1\cup F_2',\bigl(\{1,\ldots,n\}\setminus
F_1\bigr)\bigcup\bigl(\{1,\ldots,n\}\setminus F_2\bigr)'\bigr\}$.
Then, taking to account Proposition \ref{pr:Omega}, we have that
$a^{-1}[f_1]a=\{f_2\}\subseteq [f_2]$ and
$a[f_2]a^{-1}=\{f_1\}\subseteq [f_1]$, whence $\phi_{[f_1]}$ and
$\phi_{[f_2]}$ are equivalent. This completes the proof.
\end{proof}

Lemmas~\ref{lm:Ef-Tr-2}, \ref{lm:Ef-Tr-3}, \ref{lm:Ef-Tr-4},
\ref{lm:Ef-Tr-5}, \ref{lm:Ef-Tr-6}, \ref{lm:Ef-Tr-7} imply the
statement of our theorem. We are done.

\section{Definition of the ordered partition semigroup
$\IOP_n$}\label{IOP_n}

Let $n\in\mathbb{N}$. Consider the natural linear order on the set
$\{1,\ldots,n\}$. Take $A\subseteq\{1,\ldots,n\}$. Denote by $\m_A$
the minimum element of $A$ with respect to this order.

Denote by $\IOP_n$ the set of all elements $a=\bigl(A_i\cup
B_i'\bigr)_{i\in I}$ of $\IP_n$ such that
\begin{equation}\label{eq:IOP-definition}
\m_{A_i}\leq\m_{A_j}\Rightarrow\m_{B_i}\leq\m_{B_j}~\mbox{for
all}~i,j\in I.
\end{equation}
The following theorem shows that $\IOP_n$ is an inverse subsemigroup
of $\IP_n$.

\begin{theorem}\label{th:IOP_n}
$\IOP_n$ is an inverse subsemigroup of $\IP_n$.
\end{theorem}

\begin{proof}
That $a\in\IOP_n$ implies $a^{-1}\in\IOP_n$, follows immediately
from \eqref{eq:IOP-definition}. It remains to prove that $\IOP_n$ is
a subsemigroup of $\IP_n$.

Take $a,b\in\IOP_n$. Set $c=ab$. Let $a=\bigl(A_i\cup
B_i'\bigr)_{i\in I}$, $b=\bigl(C_j\cup D_j'\bigr)_{j\in J}$.
Obviously, $0\in\IOP_n$, so we may assume that $c\ne 0$. Let also
$c=\bigl(E_k\cup F_k'\bigr)_{k\in K}$ and set a linear order
$\preceq$ on $K$, given by
\begin{equation}
\m_{E_k}\leq\m_{E_l}~\mbox{if and only if}~k\preceq l~\mbox{for
all}~k,l\in K.
\end{equation}
Let now $K=\{k_1,\ldots,k_m\}$ and $k_1\preceq
k_2\preceq\ldots\preceq k_m$. Set $P_i=E_{k_i}$ and $Q_i=F_{k_i}$
for all $i$, $1\leq i\leq m$. Then we have
\begin{equation}\label{eq:wichtig}
\m_{P_1}\leq\ldots\leq\m_{P_m}.
\end{equation}
Obviously, we have that $1\equiv_{a}1'$ and $1\equiv_{b}1'$. So
$1\equiv_{c}1'$. Due to this fact, we obtain that $\{1,1'\}$ is a
subset of the block $P_1\cup Q_1'$ of the element $c$. This implies
that $\m_{Q_1}=1$. So $\m_{Q_1}\leq\m_{Q_2}$ and $\m_{Q_1}$ is the
first number among the numbers $\m_{Q_1},\ldots,\m_{Q_m}$.

Suppose now that $\m_{Q_1}\leq\ldots\leq\m_{Q_t}$ and that
$\m_{Q_1},\ldots,\m_{Q_t}$ are the first $t$ numbers among the
numbers $\m_{Q_1},\ldots,\m_{Q_m}$, for some $t$, $t<m$. Then
$\m_{Q_t}\leq\m_{Q_{t+1}}$. Since $Q_1,Q_2,\ldots,Q_t$ are all
$\lambda_{ab}$-classes, we obtain that each $Q_i$, $i\leq t$, is a
union of some $\lambda_{b}$-classes and so $Z=Q_1\cup
Q_2\cup\ldots\cup Q_t$ is a union of the sets $D_{r}$, $r\in
R\subseteq J$. Further, we have that there is a subset $U$ of $I$
such that $\bigcup\limits_{u\in U}B_u=\bigcup\limits_{r\in R}C_r=W$.
There is $r_0$ of $R$ such that $\m_{Q_t}\in D_{r_0}$. Then,
obviously, $\m_{D_{r_0}}=\m_{Q_t}$. Since $\m_{Q_1},\ldots,\m_{Q_t}$
are the first $t$ numbers among $\m_{Q_1},\ldots,\m_{Q_m}$, we
obtain that
\begin{equation}
\{1,\ldots,\m_{D_{r_0}}\}=\{1,\ldots,\m_{Q_1}\}\subseteq\bigcup\limits_{i=1}^{t}Q_i=\bigcup\limits_{r\in
R}D_r.
\end{equation}
The latter implies that $\m_{D_r}$, $r\in R$, are the first $\lmod
R\rmod$ numbers among the numbers $\m_{D_j}$, $j\in J$. Besides,
$\m_{D_{r}}\leq\m_{D_{r_0}}$ for all $r\in R$. Then, taking to
account that $b\in\IOP_n$, we obtain that
$\{1,\ldots,\m_{C_{r_0}}\}\subseteq\bigcup\limits_{r\in R}C_r$ and
$\m_{C_{r}}\leq\m_{C_{r_0}}$ for all $r\in R$. Then, taking to
account $\bigcup\limits_{u\in U}B_u=\bigcup\limits_{r\in R}C_r$, we
obtain that $\m_{B_u}$, $u\in U$, are the first $\lmod U\rmod$
numbers among the numbers $\m_{B_i}$, $i\in I$. Then, applying
$a\in\IOP_n$, we obtain that $\m_{A_u}$, $u\in U$, are the first
$\lmod U\rmod$ numbers among the numbers $\m_{A_i}$, $i\in I$. Note
that $\bigcup\limits_{u\in U}A_u=\bigcup\limits_{i=1}^{t}P_t=Y$. Put
$y=\m_{\{1,\ldots,n\}\setminus Y}$, $w=\m_{\{1,\ldots,n\}\setminus
W}$ and $y=\m_{\{1,\ldots,n\}\setminus Z}$. Then due to what we have
already obtained and due to~\eqref{eq:wichtig}, we have that
$y=\m_{P_{t+1}}$. Suppose now that $z=\m_{Q_{g}}$, $g>t$. Then due
to our assumption, we have that
\begin{multline}\label{eq:final}
\m_{Q_1},\ldots,\m_{Q_t},z~\mbox{are the first}~t+1~\mbox{numbers}\\
\mbox{among the numbers}~\m_{Q_1},\ldots,\m_{Q_m}.
\end{multline}
Due to $a,b\in\IOP_n$, we have that $y\equiv_a w'$ and $w\equiv_b
z'$, whence $y\equiv_c z'$. This implies that $z\in Q_{t+1}$, whence
$z=\m_{Q_{t+1}}$.

Thus, due to~\eqref{eq:final}, we obtain that inductive arguments
lead us to
\begin{equation}\label{eq:one-more-last}
\m_{Q_1}\leq\ldots\leq\m_{Q_m}.
\end{equation}
The conditions~\eqref{eq:wichtig} and~\eqref{eq:one-more-last}
complete the proof.
\end{proof}

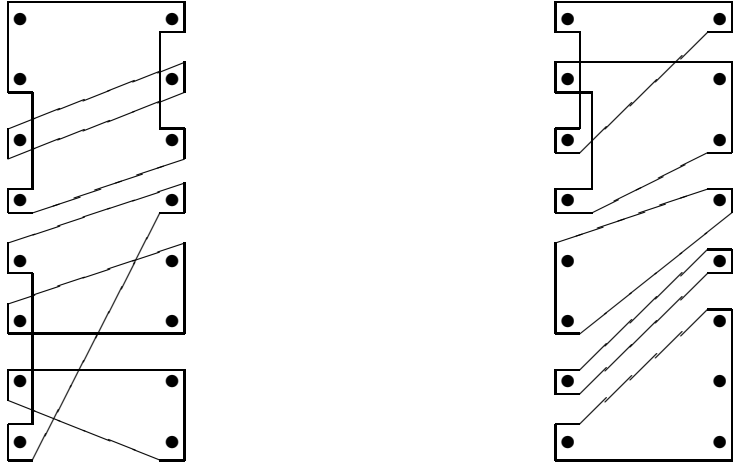
\begin{figure}
\special{em:linewidth 0.4pt} \unitlength 0.80mm
\linethickness{0.4pt}
\begin{picture}(150.00,75.00)
\put(15.00,10.00){\makebox(0,0)[cc]{$\bullet$}}
\put(15.00,20.00){\makebox(0,0)[cc]{$\bullet$}}
\put(15.00,30.00){\makebox(0,0)[cc]{$\bullet$}}
\put(15.00,40.00){\makebox(0,0)[cc]{$\bullet$}}
\put(15.00,50.00){\makebox(0,0)[cc]{$\bullet$}}
\put(15.00,60.00){\makebox(0,0)[cc]{$\bullet$}}
\put(15.00,70.00){\makebox(0,0)[cc]{$\bullet$}}
\put(15.00,00.00){\makebox(0,0)[cc]{$\bullet$}}

\put(40.00,00.00){\makebox(0,0)[cc]{$\bullet$}}
\put(40.00,10.00){\makebox(0,0)[cc]{$\bullet$}}
\put(40.00,20.00){\makebox(0,0)[cc]{$\bullet$}}
\put(40.00,30.00){\makebox(0,0)[cc]{$\bullet$}}
\put(40.00,40.00){\makebox(0,0)[cc]{$\bullet$}}
\put(40.00,50.00){\makebox(0,0)[cc]{$\bullet$}}
\put(40.00,60.00){\makebox(0,0)[cc]{$\bullet$}}
\put(40.00,70.00){\makebox(0,0)[cc]{$\bullet$}}

\put(105.00,10.00){\makebox(0,0)[cc]{$\bullet$}}
\put(105.00,20.00){\makebox(0,0)[cc]{$\bullet$}}
\put(105.00,30.00){\makebox(0,0)[cc]{$\bullet$}}
\put(105.00,40.00){\makebox(0,0)[cc]{$\bullet$}}
\put(105.00,50.00){\makebox(0,0)[cc]{$\bullet$}}
\put(105.00,60.00){\makebox(0,0)[cc]{$\bullet$}}
\put(105.00,70.00){\makebox(0,0)[cc]{$\bullet$}}
\put(105.00,00.00){\makebox(0,0)[cc]{$\bullet$}}

\put(130.00,00.00){\makebox(0,0)[cc]{$\bullet$}}
\put(130.00,10.00){\makebox(0,0)[cc]{$\bullet$}}
\put(130.00,20.00){\makebox(0,0)[cc]{$\bullet$}}
\put(130.00,30.00){\makebox(0,0)[cc]{$\bullet$}}
\put(130.00,40.00){\makebox(0,0)[cc]{$\bullet$}}
\put(130.00,50.00){\makebox(0,0)[cc]{$\bullet$}}
\put(130.00,60.00){\makebox(0,0)[cc]{$\bullet$}}
\put(130.00,70.00){\makebox(0,0)[cc]{$\bullet$}}

\drawline(13.00,-03.00)(13.00,03.00)
\drawline(13.00,03.00)(17.00,03.00)
\drawline(17.00,03.00)(17.00,28.00)
\drawline(17.00,28.00)(13.00,28.00)
\drawline(13.00,28.00)(13.00,33.00)
\drawline(13.00,33.00)(42.00,43.00)
\drawline(42.00,43.00)(42.00,38.00)
\drawline(42.00,38.00)(38.00,38.00)
\drawline(38.00,38.00)(17.00,-03.00)
\drawline(17.00,-03.00)(13.00,-03.00)

\drawline(13.00,07.00)(13.00,12.00)
\drawline(13.00,12.00)(42.00,12.00)
\drawline(42.00,12.00)(42.00,-03.00)
\drawline(42.00,-03.00)(38.00,-03.00)
\drawline(38.00,-03.00)(13.00,07.00)

\drawline(13.00,18.00)(13.00,23.00)
\drawline(13.00,23.00)(42.00,33.00)
\drawline(42.00,33.00)(42.00,18.00)
\drawline(42.00,18.00)(13.00,18.00)

\drawline(13.00,47.00)(13.00,52.00)
\drawline(13.00,52.00)(42.00,63.00)
\drawline(42.00,63.00)(42.00,58.00)
\drawline(42.00,58.00)(13.00,47.00)

\drawline(13.00,38.00)(13.00,42.00)
\drawline(13.00,42.00)(17.00,42.00)
\drawline(17.00,42.00)(17.00,58.00)
\drawline(17.00,58.00)(13.00,58.00)
\drawline(13.00,58.00)(13.00,73.00)
\drawline(13.00,73.00)(42.00,73.00)
\drawline(42.00,73.00)(42.00,68.00)
\drawline(42.00,68.00)(38.00,68.00)
\drawline(38.00,68.00)(38.00,52.00)
\drawline(38.00,52.00)(42.00,52.00)
\drawline(42.00,52.00)(42.00,47.00)
\drawline(42.00,47.00)(17.00,38.00)
\drawline(17.00,38.00)(13.00,38.00)

\drawline(103.00,-03.00)(103.00,03.00)
\drawline(103.00,03.00)(107.00,03.00)
\drawline(107.00,03.00)(128.00,22.00)
\drawline(128.00,22.00)(132.00,22.00)
\drawline(132.00,22.00)(132.00,-03.00)
\drawline(132.00,-03.00)(103.00,-03.00)

\drawline(103.00,08.00)(103.00,12.00)
\drawline(103.00,12.00)(107.00,12.00)
\drawline(107.00,12.00)(128.00,32.00)
\drawline(128.00,32.00)(132.00,32.00)
\drawline(132.00,32.00)(132.00,28.00)
\drawline(132.00,28.00)(128.00,28.00)
\drawline(128.00,28.00)(107.00,08.00)
\drawline(107.00,08.00)(103.00,08.00)

\drawline(103.00,18.00)(103.00,33.00)
\drawline(103.00,33.00)(128.00,42.00)
\drawline(128.00,42.00)(132.00,42.00)
\drawline(132.00,42.00)(132.00,38.00)
\drawline(132.00,38.00)(107.00,18.00)
\drawline(107.00,18.00)(103.00,18.00)

\drawline(103.00,48.00)(103.00,52.00)
\drawline(103.00,52.00)(107.00,52.00)
\drawline(107.00,52.00)(107.00,68.00)
\drawline(107.00,68.00)(103.00,68.00)
\drawline(103.00,68.00)(103.00,73.00)
\drawline(103.00,73.00)(132.00,73.00)
\drawline(132.00,73.00)(132.00,68.00)
\drawline(132.00,68.00)(128.00,68.00)
\drawline(128.00,68.00)(107.00,48.00)
\drawline(107.00,48.00)(103.00,48.00)

\drawline(103.00,38.00)(103.00,42.00)
\drawline(103.00,42.00)(109.00,42.00)
\drawline(109.00,42.00)(109.00,58.00)
\drawline(109.00,58.00)(103.00,58.00)
\drawline(103.00,58.00)(103.00,63.00)
\drawline(103.00,63.00)(132.00,63.00)
\drawline(132.00,63.00)(132.00,48.00)
\drawline(132.00,48.00)(128.00,48.00)
\drawline(128.00,48.00)(109.00,38.00)
\drawline(109.00,38.00)(103.00,38.00)

\end{picture}
\caption{Elements of $\IOP_8$.}\label{fig:f4}

\end{figure}

Thus, due to~Theorem \ref{th:IOP_n}, we can name $\IOP_n$ as the
\emph{inverse ordered partition semigroup} of degree $n$. On
Fig.~\ref{fig:f4} we give some examples of elements of $\IOP_8$.

Recall that a subsemigroup $T$ of a semigroup $S$ is said to be an
\emph{$\GH$-cross-section} of $S$ if $T$ contains exactly one
representative from each $\GH$-class of $S$. In the following
proposition we show that $\IOP_n$ is an $\GH$-cross-section of
$\IP_n$.

\begin{proposition}\label{pr:IOP-cross-section}
$\IOP_n$ is an $\GH$-cross-section of $\IP_n$.
\end{proposition}

\begin{proof}
Follows from \eqref{eq:IOP-definition}, Theorem~\ref{th:Green} and
Theorem~\ref{th:IOP_n}.
\end{proof}

As a consequence of Proposition~\ref{pr:IOP-cross-section}, we
obtain the following corollary.

\begin{corollary}\label{cor:IOP}
Let $n\in\mathbb{N}$. Then $E(\IOP_n)=E(\IP_n)$.
\end{corollary}

\begin{proof}
Recall that every maximal subgroup of an arbitrary semigroup $S$
coincides with some $\GH$-class of $S$, which contains an idempotent
(see~\cite{Howie}). Then every $\GH$-cross-section of $\IP_n$
contains all the idempotents of $\IP_n$. In particular,
$E(\IOP_n)=E(\IP_n)$, which was required.
\end{proof}

\section{Acknowledgments}

The author is indebted to Professor Norman Reilly and to the two
anonymous referees whose comments and suggestions contributed to a
significant improvement of this paper.

\noindent Mathematical Institute, University of St Andrews\\
St Andrews KY16 9SS, Scotland\\
e-mail: {\tt victor\symbol{64}mcs.st-and.ac.uk}


\begin{thebibliography}{0}

\bibitem{AO} M. Aguiar, R. C. Orellana, The Hopf Algebra of Uniform Block
Permutations, 17th International Conference on Formal Power Series
and Algebraic Combinatorics, Taormina , July 2005.

\bibitem{Brauer} R. Brauer, On Algebras Which are Connected with the
Semisimple Continuous Groups, {\it Ann Math.} {\bf 38}(2) (1937)
857-872.

\bibitem{Changchang} Xi. Changchang, Partition algebras are cellular,
{\it Compositio Math.} {\bf 119} (1999) 99-109.

\bibitem{Cowan-Reilly} D. Cowan, N. Reilly, Partial cross-sections
of symmetric inverse semigroups, {\it Internat. J. Algebra Comput.}
{\bf 5}(3) (1995) 259-287.

\bibitem{F} D. G. FitzGerald, A presentation for the monoid of uniform block permutations,
{\it Bull. Austral. Math. Soc.} {\bf 68} (2003) 317-324.

\bibitem{FL} D. G. FitzGerald, J. Leech, Dual symmetric inverse
monoids and representation theory, {\it J. Austral. Math.
Soc.(Series A)} {\bf 64} (1998) 345-367.

\bibitem{Howie} J. M. Howie, {\it Fundamentals of
Semigroup Theory} (Oxford/Clarendon Press, 1995).

\bibitem{Kerov} S. Kerov, Realizations of representations of the Brauer semigroup,
{\it Zap. Nauchn. Sem. LOMI} {\bf 164} (1987) 188-193.

\bibitem{KM} G. Kudryavtseva, V. Maltcev, On the structure of two generalizations
of the full symmetric inverse semigroup, Preprint, Kyiv University.

\bibitem{KMM} G. Kudryavtseva, V. Maltcev, V. Mazorchuk, $\GL$- and $\GR$-cross-sections
in the Brauer semigroup, {\it Semigroup Forum} {\bf 72} (2006)
223-248.

\bibitem{Kurosh} A. G. Kurosh, {\it Group Theory} (Moscow/Science, 1967) (in Russian).

\bibitem{LF} K. W. Lau, D. G. FitzGerald, Ideal structure of the Kauffman and related
monoids, to appear in {\it Comm. Algebra}.

\bibitem{Lawson} M. V. Lawson, {\it Inverse Semigroups: The Theory of
Partial Symmetries} (Singapore/World Scientific, 1998).

\bibitem{Maltcev} V. Maltcev, Ideals and systems of generators in the Brauer semigroup
$\mathfrak{B}_n$, {\it Reports of Kyiv University} {\bf 2}(2) (2004)
59-66.

\bibitem{Maltcev-PB-C} V. Maltcev, Cross-sections of Green relations and retracts of
semigroups $\mathcal{P}\mathfrak{B}\sb n$ and $\mathfrak{C}\sb n$,
{\it Scientific Proceedings of the Kyiv-Mohyla Academy} {\bf 39}
(2005) 11-24.

\bibitem{Maltcev-IP} V. Maltcev, On one inverse subsemigroup of
$\mathfrak{C}_n$, to appear.

\bibitem{MM} V. Maltcev, V. Mazorchuk, Presentation of the singular part
of the Brauer monoid, accepted for publication in {\it Mathematicae
Bogemica}.

\bibitem{Mazorchuk-PB} V. Mazorchuk, On the structure of the Brauer semigroup and its
partial analogue, {\it Problems in Algebra} {\bf 13} (1998) 29-45.

\bibitem{Mazorchuk-End} V. Mazorchuk, Endomorphisms of $\mathfrak{B}\sb n,
\mathcal{P}\mathfrak{B}\sb n$, and $\mathfrak{C}\sb n$, {\it
Communication in Algebra} {\bf 30}(7) (2002) 3489-3513.

\bibitem{Meakin} J. Meakin, On the structure of inverse semigroups, {\it
Semigroup Forum} {\bf 12} (1976) 6-14.

\bibitem{Petrich} M. Petrich, {\it Inverse semigroups} (New York/Wiley \& Sons, 1984).

\end{thebibliography}
\end{document}